\RequirePackage{fix-cm}
\documentclass[smallextended]{svjour3}
\smartqed 
\usepackage{graphicx}
\usepackage{graphicx}
\usepackage{amsmath}
\usepackage{amssymb}
\usepackage{latexsym}
\usepackage{mathrsfs}
\usepackage{colortbl}
\usepackage{amscd}
\usepackage{tikz-cd}
\usepackage{comment}

\def\diam{\mathop{\mathrm{diam}}}
\def\diag{\mathop{\mathrm{diag}}}

\def\div{\mathop{\mathrm{div}}}
\def\divh{\mathop{\mathrm{div}_h}}
\def\curl{\mathop{\mathrm{curl}}}
\def\curlh{\mathop{\mathrm{curl}_h}}
\def\rot{\mathop{\mathrm{rot}}}

\def\<{\mathop{\textless}}
\def\>{\mathop{\textgreater}}

\def\card{\mathop{\rm{card}}}

\def\Int{\mathop{\rm{int}}}

\spnewtheorem{thr}{Theorem}{\bf}{\it}
\spnewtheorem{defi}{Definition}{\bf}{\it}
\spnewtheorem{lem}{Lemma}{\bf}{\it}
\spnewtheorem{coro}{Corollary}{\bf}{\it}
\spnewtheorem{assume}{Assumption}{\bf}{\it}
\spnewtheorem{ex}{Example}{\bf}{\it}
\spnewtheorem{Case}{Case}{\bf}{\it}
\spnewtheorem*{pf*}{Proof}{\bf}{\rm}
\spnewtheorem*{rem*}{Remark:}{\it}{\it}
\spnewtheorem*{ex*}{Example:}{\it}{\it}
\spnewtheorem{Cond}{Condition}{\bf}{\it}
\spnewtheorem{rem}{Remark}{\it}{\it}

\spnewtheorem*{lem1*}{Lemma 1}{\bf}{\rm}
\spnewtheorem*{lem3*}{Lemma 3}{\bf}{\rm}

\spnewtheorem*{th21*}{Theorem A}{\bf}{\it}
\spnewtheorem*{thrBA*}{Theorem B}{\bf}{\it}
\spnewtheorem*{thrF*}{Theorem D}{\bf}{\it}
\spnewtheorem*{thrin*}{Theorem C}{\bf}{\it}
\spnewtheorem*{thrlp*}{Theorem E}{\bf}{\it}

\spnewtheorem*{th2*}{Theorem 2}{\bf}{\it}

\spnewtheorem*{coex1*}{Counterexample 1}{\it}{\it}
\spnewtheorem*{coex2*}{Counterexample 2}{\it}{\it}

\allowdisplaybreaks[3]

\newcounter{sone}
\newcounter{stwo}
\newcounter{sthree}
\newcounter{sfour}
\newcounter{sfive}
\newcounter{ssix}
\newcounter{lone}
\newcounter{ltwo}
\newcounter{lthree}
\newcounter{lfour}
\newcounter{lfive}
\newcounter{lsix}
\makeatletter
    
    \@addtoreset{equation}{section}
  \makeatother

\begin{document}

\title{Morley finite element analysis for fourth-order elliptic equations {under a semi-regular mesh condition}
}

\titlerunning{Morley FE analysis for fourth-order elliptic equations}  

\author{Hiroki Ishizaka
}

\institute{Hiroki Ishizaka \at
                Team FEM, Matsuyama, Japan \\
              \email{h.ishizaka005@gmail.com}
}

\date{Received: date / Accepted: date}

\maketitle

\begin{abstract}

In this study, we present a precise anisotropic interpolation error estimate for the Morley finite element method (FEM) and apply it to fourth-order elliptic equations. We do not impose the shape-regularity mesh condition in the analysis. Anisotropic meshes can be used for this purpose. The main contributions of this study include providing a new proof of the term consistency. This enables us to obtain an anisotropic consistency error estimate. The core idea of the proof involves using the relationship between the Raviart--Thomas and Morley finite-element spaces. Our results indicate optimal convergence rates and imply that the modified Morley FEM may be effective for errors.

\end{abstract}

\section{Introduction}
This study investigates the error estimates of nonconforming finite element methods (FEMs) for fourth-order {elliptic} equations {under a semi-regular mesh condition. }

In the context of standard FEMs, the shape-regular family of triangulations, in which triangles or tetrahedra cannot be overly flat, is widely used to estimate the optimal order errors. Moreover, anisotropic meshes may be effective for problems in which the solution exhibits anisotropic behaviour in certain domain directions, such as singularly perturbed differential equations with boundary or interior layers, differential equations with edge singularities, and flow problems such as the Stokes and Navier--Stokes equations. For these problems, using a regular mesh may contribute to an increase in errors. However, anisotropic meshes do not satisfy the shape-regularity condition. Therefore, it is necessary to extend the previous theory of error analysis. Several anisotropic element methods have recently been developed \cite{Ape99,ApeDob92,Ish22b,IshKobTsu23a}. {These methods aim to obtain optimal error estimates under the \textit{semi-regular condition} defined in Assumption \ref{ass1}, see \cite{IshKobTsu23a}, or the \textit{maximum-angle condition}, which allows the use of anisotropic meshes, see \cite{BabAzi76} for two-dimensional and \cite{Kri92} for three-dimensional cases. }

In this study, we consider an anisotropic Morley FEM for fourth-order {elliptic} equations. When we consider fourth-order {elliptic} problems, the weak solution belongs to the Sobolev space $H^2$, which implies that finite element spaces must be used to approximate the $H^2$ space. However, in conforming cases, we must use polynomials of piecewise degree five or higher (e.g.  \cite[p. 334]{Cia02}) or the Hsieh--Clough--Tocher method (e.g. \cite[p. 340]{Cia02}), which is a macroelement technique. To overcome this difficulty, the Morley FEM {is} considered owing to its low degrees of freedom. Because the Morley finite element space belongs to neither $H^2$ nor $H^1$, an error analysis is difficult to conduct (e.g. \cite{ArnBre85,LasLea75,Mol68,NilTaiWin01,Ran79a,Ran79b}).

Studies on anisotropic Morley FEMs that do not impose the shape-regularity condition for the mesh partition include \cite{MaoNicShi10}. A two-dimensional case was considered numerically for a problem with boundary layers. {In this paper, we present} a precise anisotropic interpolation error estimate of an alternative approach for Morley FEM in the broken $H^2$ norm. Furthermore, we present a Crouzeix--Raviart (CR) finite-element interpolation error estimate for anisotropic meshes. Obtaining the CR interpolation error estimate is not a novel concept. However, it {is} introduced for comparison with the Morley interpolation error. The Morley finite element method (FEM) is nonconforming. Therefore, the error between the exact and approximate Morley finite element solutions with an energy norm {is} divided into two parts. One is the optimal approximation error in the Morley finite element space, and the other is the consistency error term. The former involves Morley interpolation errors (Theorem \ref{thr2}). However, estimating the consistency error term for anisotropic meshes is difficult. The standard argument uses scaling arguments and the trace theorem. In particular, trace inequality on anisotropic meshes does not lead to optimal order error estimates. We leverage the relationship between the Raviart--Thomas (RT) and Morley finite element spaces to overcome this difficulty. We {consider} the usual and modified Morley finite-element spaces to apply fourth-order {elliptic} problems and stream function formulations. We theoretically {compare} these errors. {However, we impose the regularity assumption that is stronger than the original Arnold--Brezzi modified Morley method \cite{ArnBre85}}.

We have organized the remainder of this paper as follows. In Section 2, we introduce the basic notation. Section 3 introduces the anisotropic $L^2$-orthogonal projection error estimate. Section 4 introduces nonconforming error estimates. Section 5 presents useful relations and existing interpolation error estimates for the analysis. Section 6 presents the modified Morley finite element method for a fourth-order {elliptic} problem, and Section 7 introduces the stream function formulation. {The applications in Sections 6 and 7 have been studied in a two-dimensional case, where difficulties in the three-dimensional case are described in Section 6.6.}

In this paper, we {use} standard Sobolev spaces with associated norms (e.g., see \cite{ErnGue04,GirRav86,Gri11}). Let $T \subset \mathbb{R}^d$, $d \in \{1,2,3\}$ be a simplex. For $k \in \mathbb{N}_0:= \mathbb{N} \cup \{ 0 \}$, $\mathbb{P}^k(T)$ is spanned by the restriction to $T$ of polynomials in $\mathbb{P}^k$, where $\mathbb{P}^k$ denotes the space of polynomials with a degree of at most $k$. We set $N^{(d,k)} := \dim (\mathbb{P}^k) = \begin{pmatrix}
 d+k \\
 k
\end{pmatrix}.
$
Throughout, we denote by $c$ a constant independent of $h$ (defined later) and of the angles and aspect ratios of simplices unless specified otherwise, and all constants $c$ are bounded {from above} if the maximum angle is bounded {from above}. These values vary across different contexts.

\section{Preliminaries} \label{Pre=Sec}
We now introduce a Jensen-type inequality (see \cite[Exercise 12.1]{ErnGue21a}). Let $r,s$ be two nonnegative real numbers and $\{ x_i \}_{i \in I}$ be a finite sequence of nonnegative numbers. It then holds that
\begin{align}
\displaystyle
\begin{cases}
\left( \sum_{i \in I} x_i^s \right)^{\frac{1}{s}} \leq \left( \sum_{i \in I} x_i^r \right)^{\frac{1}{r}} \quad \text{if $r \leq s$},\\
\left( \sum_{i \in I} x_i^s \right)^{\frac{1}{s}} \leq \card(I)^{\frac{r-s}{rs}} \left( \sum_{i \in I} x_i^r \right)^{\frac{1}{r}} \quad \text{if $r \> s$}.
\end{cases} \label{jensen}
\end{align}

\textbf{Meshes, mesh faces, and jumps.} Let $\Omega \subset \mathbb{R}^d${, $d \in \{ 2,3\}$,} be a bounded polyhedral domain, and $\mathbb{T}_h = \{ T \}$ be a simplicial mesh of $\overline{\Omega}$ comprising closed {$d$-simplices} such that
\begin{align*}
\displaystyle
\overline{\Omega} = \bigcup_{T \in \mathbb{T}_h} T,
\end{align*}
where $h := \max_{T \in \mathbb{T}_h} h_{T}$ and $ h_{T} := \diam(T)$. For simplicity, we assume that $\mathbb{T}_h$ is conformal; that is, $\mathbb{T}_h$ is a simplicial mesh of $\overline{\Omega}$ without hanging nodes. We define the broken (piecewise) Sobolev space as
\begin{align*}
\displaystyle
H^m(\mathbb{T}_h) &:= \left\{ \varphi \in L^2(\Omega); \ \varphi|_{T} \in H^m(T) \ \forall T \in \mathbb{T}_h  \right\}, \quad m \in \mathbb{N}_0
\end{align*}
with the norm
\begin{align*}
\displaystyle
| \varphi |_{H^m(\mathbb{T}_h)} &:= \left( \sum_{T \in \mathbb{T}_h} | \varphi |^2_{H^m(T)} \right)^{1/2} \quad \varphi \in H^m(\mathbb{T}_h).
\end{align*}

Let $Ne$ be the number of elements included in mesh $\mathbb{T}_h$. Thus, we have $\mathbb{T}_h = \{ T_j\}_{j=1}^{Ne}$. Let $\mathcal{F}_h^i$ be the set of interior faces, and $\mathcal{F}_h^{\partial}$ be the set of faces on boundary $\partial \Omega$. We set $\mathcal{F}_h := \mathcal{F}_h^i \cup \mathcal{F}_h^{\partial}$. For any $F \in \mathcal{F}_h$, we define the unit normal $n_F$ to $F$ as follows: (\roman{sone}) If  $F \in \mathcal{F}_h^i$ with $F = T_{\natural} \cap T_{\sharp}$, $T_{\natural},T_{\sharp} \in \mathbb{T}_h$, $\natural > \sharp$, let $n_F$ be the unit normal vector from $T_{\natural}$ to  $ T_{\sharp}$.  (\roman{stwo}) If $F \in \mathcal{F}_h^{\partial}$, $n_F$ is {the unit outward normal} $n$ to $\partial \Omega$. 

Let $\varphi \in H^1(\mathbb{T}_h)$. Suppose that $F \in \mathcal{F}_h^i$ with $F = T_{\natural} \cap T_{\sharp}$, $T_{\natural},T_{\sharp} \in \mathbb{T}_h$, $\natural > \sharp$. We set $\varphi_{\natural} := \varphi{|_{T_{\natural}}}$ and $\varphi_{\sharp} := \varphi{|_{T_{\sharp}}}$. The jump in $\varphi$ across $F$ is defined as
\begin{align*}
\displaystyle
[\! [ \varphi ]\!] := [\! [ \varphi ]\!]_F := \varphi_{\natural} - \varphi_{\sharp}, \quad \natural > \sharp.
\end{align*}
For a boundary face $F \in \mathcal{F}_h^{\partial}$ with $F = \partial T \cap \partial \Omega$, $[\![\varphi ]\!]_F := \varphi|_{T}$. For any $v \in H^1(\mathbb{T}_h)^d$, the notation
\begin{align*}
\displaystyle
[\![ v \cdot n ]\!] := [\![ v \cdot n ]\!]_F := v_{\natural} \cdot n_F - v_{\sharp} \cdot n_F, \quad \natural > \sharp
\end{align*}
denotes the jump in the normal component of $v$. 

We define a broken $H(\div;T)$ space as follows:
\begin{align*}
\displaystyle
H(\div;\mathbb{T}_h) := \left \{ v \in L^2(\Omega)^d; \ v |_{T} \in H(\div;T) \ \forall T \in \mathbb{T}_h  \right\}.
\end{align*}
Thus, the broken divergence operator $\divh : H(\div;\mathbb{T}_h) \to L^2(\Omega)$ is such that, for all $v \in H(\div;\mathbb{T}_h)$,
\begin{align*}
\displaystyle
(\divh v)|_{T} := \div (v |_{T}) \quad \forall T \in \mathbb{T}_h.
\end{align*}

\section{$L^2$-orthogonal projection error estimate}

\subsection{Reference elements} \label{reference}
We now define the reference elements $\widehat{T} \subset \mathbb{R}^d$.

\subsubsection*{Two-dimensional case} \label{reference2d}
Let $\widehat{T} \subset \mathbb{R}^2$ be a reference triangle with vertices $\hat{p}_1 := (0,0){^{\top}}$, $\hat{p}_2 := (1,0){^{\top}}$, and $\hat{p}_3 := (0,1){^{\top}}$. 

\subsubsection*{Three-dimensional case} \label{reference3d}
In the three-dimensional case, we consider the following two cases: (\roman{sone}) and (\roman{stwo}); see Condition \ref{cond2}.

Let $\widehat{T}_1$ and $\widehat{T}_2$ be reference tetrahedra with the following vertices:
\begin{description}
   \item[(\roman{sone})] $\widehat{T}_1$ has vertices $\hat{p}_1 := (0,0,0){^{\top}}$, $\hat{p}_2 := (1,0,0){^{\top}}$, $\hat{p}_3 := (0,1,0){^{\top}}$, and $\hat{p}_4 := (0,0,1)^T$;
 \item[(\roman{stwo})] $\widehat{T}_2$ has vertices $\hat{p}_1 := (0,0,0){^{\top}}$, $\hat{p}_2 := (1,0,0){^{\top}}$, $\hat{p}_3 := (1,1,0){^{\top}}$, and $\hat{p}_4 := (0,0,1)^T$.
\end{description}
Therefore, we set $\widehat{T} \in \{ \widehat{T}_1 , \widehat{T}_2 \}$. 
Note that the case (\roman{sone}) is called \textit{the regular vertex property}, see \cite{AcoDur99}.

\subsection{Affine mappings} \label{Affinedef}
We {introduce} a new strategy in \cite[Section 2]{IshKobTsu23a} to use anisotropic mesh partitions. {To an affine simplex $T \subset \mathbb{R}^d$, we construct two affine mappings $\Phi_{\widetilde{T}}: \widehat{T} \to \widetilde{T}$ and $\Phi_{T}: \widetilde{T} \to T$.} First, we define the affine mapping $\Phi_{\widetilde{T}}: \widehat{T} \to \widetilde{T}$ as
\begin{align}
\displaystyle
\Phi_{\widetilde{T}}: \widehat{T} \ni \hat{x} \mapsto \tilde{x} := \Phi_{\widetilde{T}}(\hat{x}) := {A}_{\widetilde{T}} \hat{x} \in  \widetilde{T}, \label{aff=1}
\end{align}
where ${A}_{\widetilde{T}} \in \mathbb{R}^{d \times d}$ is an invertible matrix. Section \ref{sec221} provides the details. We then define the affine mapping $\Phi_{T}: \widetilde{T} \to T$  as follows:
\begin{align}
\displaystyle
\Phi_{T}: \widetilde{T} \ni \tilde{x} \mapsto x := \Phi_{T}(\tilde{x}) := {A}_{T} \tilde{x} + b_{T} \in T, \label{aff=2}
\end{align}
where $b_{T} \in \mathbb{R}^d$ is a vector and ${A}_{T} \in O(d)$ is {a} rotation and mirror imaging matrix. Section \ref{sec322} provides the details. We define the affine mapping $\Phi: \widehat{T} \to T$ as
\begin{align*}
\displaystyle
\Phi := {\Phi}_{T} \circ {\Phi}_{\widetilde{T}}: \widehat{T} \ni \hat{x} \mapsto x := \Phi (\hat{x}) =  ({\Phi}_{T} \circ {\Phi}_{\widetilde{T}})(\hat{x}) = {A} \hat{x} + b_{T} \in T, 
\end{align*}
where ${A} := {A}_{T} {A}_{\widetilde{T}} \in \mathbb{R}^{d \times d}$.

\subsubsection{Construct mapping $\Phi_{\widetilde{T}}: \widehat{T} \to \widetilde{T}$} \label{sec221}
We consider affine mapping \eqref{aff=1}. We define the matrix $ {A}_{\widetilde{T}} \in \mathbb{R}^{d \times d}$ as follows: We first define the diagonal matrix as
\begin{align}
\displaystyle
\widehat{A} :=  \diag (h_1,\ldots,h_d), \quad h_i \in \mathbb{R}_+ \quad \forall i,\label{aff=3}
\end{align}
where $\mathbb{R}_+$ denotes the set of positive real numbers.

For $d=2$, we define the regular matrix $\widetilde{A} \in \mathbb{R}^{2 \times 2}$ as:
\begin{align}
\displaystyle
\widetilde{A} :=
\begin{pmatrix}
1 & s \\
0 & t \\
\end{pmatrix}, \label{aff=4}
\end{align}
with parameters
\begin{align*}
\displaystyle
s^2 + t^2 = 1, \quad t \> 0.
\end{align*}
For reference element $\widehat{T}$, let $\mathfrak{T}^{(2)}$ be {the family of triangles}
\begin{align*}
\displaystyle
\widetilde{T} &= \Phi_{\widetilde{T}}(\widehat{T}) = {A}_{\widetilde{T}} (\widehat{T}), \quad {A}_{\widetilde{T}} := \widetilde {A} \widehat{A}
\end{align*}
with vertices $\tilde{p}_1 := (0,0){^{\top}}$, $\tilde{p}_2 := (h_1,0){^{\top}}$, and $\tilde{p}_3 :=(h_2 s , h_2 t){^{\top}}$. Then,  $h_1 = |\tilde{p}_1 - \tilde{p}_2| \> 0$ and $h_2 = |\tilde{p}_1 - \tilde{p}_3| \> 0$. 

For $d=3$, we define the regular matrices $\widetilde{A}_1, \widetilde{A}_2 \in \mathbb{R}^{3 \times 3}$ as 
\begin{align}
\displaystyle
\widetilde{A}_1 :=
\begin{pmatrix}
1 & s_1 & s_{21} \\
0 & t_1  & s_{22}\\
0 & 0  & t_2\\
\end{pmatrix}, \
\widetilde{A}_2 :=
\begin{pmatrix}
1 & - s_1 & s_{21} \\
0 & t_1  & s_{22}\\
0 & 0  & t_2\\
\end{pmatrix} \label{aff=5}
\end{align}
with parameters
\begin{align*}
\displaystyle
\begin{cases}
s_1^2 + t_1^2 = 1, \ s_1 \> 0, \ t_1 \> 0, \ h_2 s_1 \leq h_1 / 2, \\
s_{21}^2 + s_{22}^2 + t_2^2 = 1, \ t_2 \> 0, \ h_3 s_{21} \leq h_1 / 2.
\end{cases}
\end{align*}
Therefore, we set $\widetilde{A} \in \{ \widetilde{A}_1 , \widetilde{A}_2 \}$. For the reference elements $\widehat{T}_i$, $i=1,2$, let $\mathfrak{T}_i^{(3)}$, $i=1,2$ be {the family of tetrahedra}
\begin{align*}
\displaystyle
\widetilde{T}_i &= \Phi_{\widetilde{T}_i} (\widehat{T}_i) =  {A}_{\widetilde{T}_i} (\widehat{T}_i), \quad {A}_{\widetilde{T}_i} := \widetilde {A}_i \widehat{A}, \quad i=1,2,
\end{align*}
with vertices
\begin{align*}
\displaystyle
&\tilde{p}_1 := (0,0,0){^{\top}}, \ \tilde{p}_2 := (h_1,0,0){^{\top}}, \ \tilde{p}_4 := (h_3 s_{21}, h_3 s_{22}, h_3 t_2){^{\top}}, \\
&\begin{cases}
\tilde{p}_3 := (h_2 s_1 , h_2 t_1 , 0){^{\top}} \quad \text{for case (\roman{sone})}, \\
\tilde{p}_3 := (h_1 - h_2 s_1, h_2 t_1,0){^{\top}} \quad \text{for case (\roman{stwo})}.
\end{cases}
\end{align*}
Subsequently, $h_1 = |\tilde{p}_1 - \tilde{p}_2| \> 0$, $h_3 = |\tilde{p}_1 - \tilde{p}_4| \> 0$, and
\begin{align*}
\displaystyle
h_2 =
\begin{cases}
|\tilde{p}_1 - \tilde{p}_3| \> 0  \quad \text{for case (\roman{sone})}, \\
|\tilde{p}_2 - \tilde{p}_3| \> 0  \quad \text{for case (\roman{stwo})}.
\end{cases}
\end{align*}

\subsubsection{Construct mapping $\Phi_{T}: \widetilde{T} \to T$}  \label{sec322}
We determine the affine mapping \eqref{aff=2} as follows: Let ${T} \in \mathbb{T}_h$ have vertices ${p}_i$ ($i=1,\ldots,d+1$). Let $b_{T} \in \mathbb{R}^d$ be the vector and ${A}_{T} \in O(d)$ be the rotation and mirror imaging matrix such that
\begin{align*}
\displaystyle
p_{i} = \Phi_T (\tilde{p}_i) = {A}_{T} \tilde{p}_i + b_T, \quad i \in \{1, \ldots,d+1 \},
\end{align*}
where the vertices $p_{i}$ ($i=1,\ldots,d+1$) satisfy the following conditions:

\begin{Cond}[Case in which $d=2$] \label{cond1}
We assume that $\overline{{p}_2 {p}_3}$ is the longest edge of ${T}$; that is, $ h_{{T}} := |{p}_2 - {p}_ 3|$. We {assume} that $h_2 \leq h_1$. We then have $h_1 = |{p}_1 - {p}_2|$ and $h_2 = |{p}_1 - {p}_3|$. {Because $\frac{1}{2} h_T < h_1 \leq h_T$, ${h_1 \approx h_T}$.}
\end{Cond}

\begin{Cond}[Case in which $d=3$] \label{cond2}
Let ${L}_i$ ($1 \leq i \leq 6$) be an edge of ${T}$. ${L}_{\min}$ denotes the edge of ${T}$ with minimum length, that is, $|{L}_{\min}| = \min_{1 \leq i \leq 6} |{L}_i|$. We set $h_2 := |{L}_{\min}|$ and assume that 
\begin{align*}
\displaystyle
&\text{the endpoints of ${L}_{\min}$ are either $\{ {p}_1 , {p}_3\}$ or $\{ {p}_2 , {p}_3\}$}.
\end{align*}
Among the four edges sharing an endpoint with ${L}_{\min}$, we consider the longest edge, ${L}^{({\min})}_{\max}$. Let ${p}_1$ and ${p}_2$ be the endpoints of edge ${L}^{({\min})}_{\max}$. Hence, we have 
\begin{align*}
\displaystyle
h_1 = |{L}^{(\min)}_{\max}| = |{p}_1 - {p}_2|.
\end{align*}
We consider cutting $\mathbb{R}^3$ with a plane that contains the midpoint of the edge ${L}^{(\min)}_{\max}$ and is perpendicular to the vector ${p}_1 - {p}_2$. We then have two cases: 
\begin{description}
  \item[(Type \roman{sone})] ${p}_3$ and ${p}_4$  belong to the same half-space;
  \item[(Type \roman{stwo})] ${p}_3$ and ${p}_4$  belong to different half-spaces.
\end{description}
In each case, we set
\begin{description}
  \item[(Type \roman{sone})] ${p}_1$ and ${p}_3$ as the endpoints of ${L}_{\min}$, that is, $h_2 =  |{p}_1 - {p}_3| $;
  \item[(Type \roman{stwo})] ${p}_2$ and ${p}_3$ as the endpoints of ${L}_{\min}$, that is, $h_2 =  |{p}_2 - {p}_3| $.
\end{description}
Finally, we have $h_3 = |{p}_1 - {p}_4|$. We implicitly assume that ${p}_1$ and ${p}_4$ belong to the same half-space. 
\end{Cond}

\begin{note}
As an example, we {define} the matrices $A_{T}$ as 
\begin{align*}
\displaystyle
A_{T} := 
\begin{pmatrix}
\cos \theta  & - \sin \theta \\
 \sin \theta & \cos \theta
\end{pmatrix}, \quad 
{A}_{T} := 
\begin{pmatrix}
 \cos \theta  & - \sin \theta & 0\\
 \sin \theta & \cos \theta & 0 \\
 0 & 0 & 1 \\
\end{pmatrix},
\end{align*}
where $\theta$ denotes the angle. 
\end{note}

\begin{note}
None of the lengths of the edges of a simplex or the measures of the simplex {are} changed by the transformation. 	
\end{note}

\subsection{Additional notation and assumption} \label{addinot}
\begin{figure}[tbhp]
\vspace{-7cm}
  \includegraphics[bb=0 0 944 702,scale=0.55]{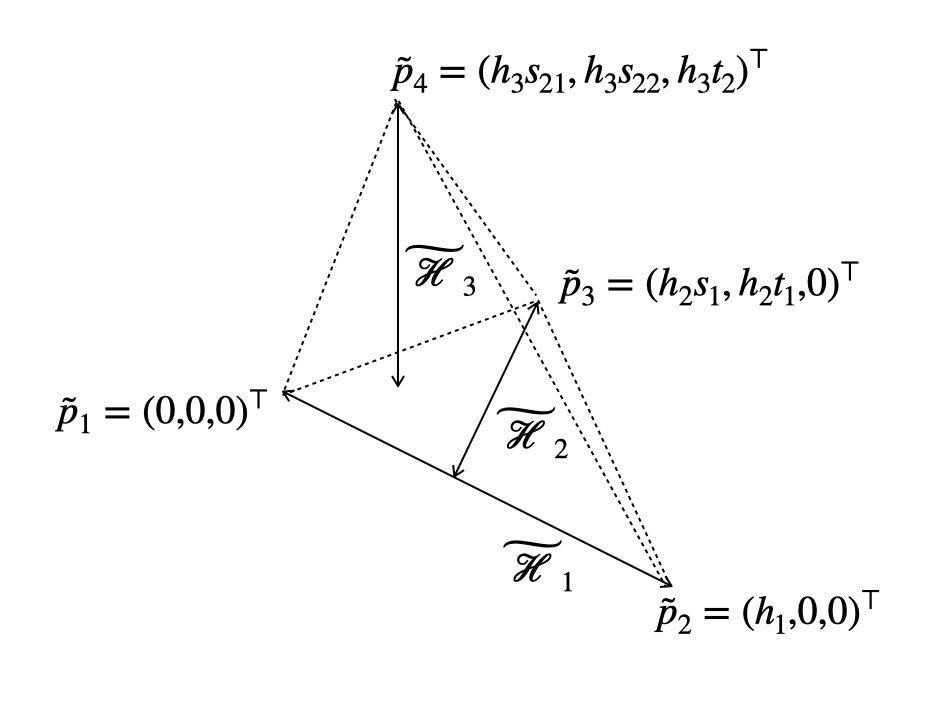}
\caption{New parameters $\widetilde{\mathscr{H}}_i$, $i=1,2,3$}
\label{mathscrH}
\end{figure}

For convenience, we {introduce} the following additional notation.
We {define} a parameter $\widetilde{\mathscr{H}}_i$, $i=1,\ldots,d$, as
\begin{align*}
\displaystyle
\begin{cases}
\widetilde{\mathscr{H}}_1 := h_1, \quad \widetilde{\mathscr{H}}_2 := h_2 t \quad \text{if $d=2$}, \\
\widetilde{\mathscr{H}}_1 := h_1, \quad \widetilde{\mathscr{H}}_2 := h_2 t_1, \quad \widetilde{\mathscr{H}}_3 := h_3 t_2 \quad \text{if $d=3$},
\end{cases}
\end{align*}
{see Fig. \ref{mathscrH}.}

\begin{assume} \label{assume0}
In an anisotropic interpolation error analysis, we {impose} a geometric condition for the simplex $\widetilde{T}$:
\begin{enumerate}
 \item If $d=2$, there are no additional conditions;
 \item If $d=3$, there exists a positive constant $M$ independent of $h_{\widetilde{T}}$ such that $|s_{22}| \leq M \frac{h_2 t_1}{h_3}$. Note that if $s_{22} \neq 0$, this condition means that the order concerning $h_T$ of $h_3$ coincides with the order of $h_2$, and if $s_{22} = 0$, the order of $h_3$ may be different from that of $h_2$. 
\end{enumerate}
\end{assume}

We {define} the vectors ${r}_n \in \mathbb{R}^d$ and $n=1,\ldots,d$ as follows: If $d=2$,
\begin{align*}
\displaystyle
{r}_1 := \frac{p_2 - p_1}{|p_2 - p_1|}, \quad \tilde{r}_2 := \frac{p_3 - p_1}{|p_3 - p_1|},
\end{align*}
and if $d=3$,
\begin{align*}
\displaystyle
&{r}_1 := \frac{p_2 - p_1}{|p_2 - p_1|}, \quad {r}_3 := \frac{p_4 - p_1}{|p_4 - p_1|}, \quad
\begin{cases}
\displaystyle
{r}_2 := \frac{p_3 - p_1}{|p_3 - p_1|}, \quad \text{for case (\roman{sone})}, \\
\displaystyle
{r}_2 := \frac{p_3 - p_2}{|p_3 - p_2|} \quad \text{for case (\roman{stwo})}.
\end{cases}
\end{align*}
Furthermore, we define the vectors $\tilde{r}_n \in \mathbb{R}^d$ and $n=1,\ldots,d$ as follows. If $d=2$,
\begin{align*}
\displaystyle
\tilde{r}_1 := (1 , 0){^{\top}}, \quad \tilde{r}_2 := (s,t){^{\top}},
\end{align*}
and if $d=3$,
\begin{align*}
\displaystyle
&\tilde{r}_1 := (1 , 0,0){^{\top}}, \quad \tilde{r}_3 := ( s_{21}, s_{22} , t_2){^{\top}}, \quad
\begin{cases}
\tilde{r}_2 := ( s_1 ,  t_1 , 0){^{\top}} \quad \text{for case (\roman{sone})}, \\
\tilde{r}_2 := (- s_1,  t_1,0){^{\top}} \quad \text{for case (\roman{stwo})}.
\end{cases}
\end{align*}

For a sufficiently smooth function $\varphi$ and vector function $v := (v_{1},\ldots,v_{d}){^{\top}}$, we define the directional derivative of $i \in \{ 1, \ldots, d \}$ as:
\begin{align*}
\displaystyle
\frac{\partial \varphi}{\partial {r_i}} &:= ( {r}_i \cdot  \nabla_{x} ) \varphi = \sum_{i_0=1}^d ({r}_i)_{i_0} \frac{\partial \varphi}{\partial x_{i_0}^{}}, \\
\frac{\partial v}{\partial r_i} &:= \left(\frac{\partial v_{1}}{\partial r_i}, \ldots, \frac{\partial v_{d}}{\partial r_i} \right){^{\top}} 
= ( ({r}_i  \cdot \nabla_{x}) v_{1}, \ldots, ({r}_i  \cdot \nabla_{x} ) v_{d} ){^{\top}}.
\end{align*}
For a multi-index $\beta = (\beta_1,\ldots,\beta_d) \in \mathbb{N}_0^d$, we use the following notation.
\begin{align*}
\displaystyle
\partial^{\beta} \varphi := \frac{\partial^{|\beta|} \varphi}{\partial x_1^{\beta_1} \ldots \partial x_d^{\beta_d}}, \quad \partial^{\beta}_{r} \varphi := \frac{\partial^{|\beta|} \varphi}{\partial r_1^{\beta_1} \ldots \partial r_d^{\beta_d}}.
\end{align*}
Note that $\partial^{\beta} \varphi \neq  \partial^{\beta}_{r} \varphi$.

We proposed a new geometric parameter $H_{T}$ in \cite{IshKobTsu21a}.
 \begin{defi} \label{defi1}
 Parameter $H_T$ is defined as follows:
\begin{align*}
\displaystyle
H_T := \frac{\prod_{i=1}^d h_i}{|T|} h_T,
\end{align*}
{where $|{T}|$ is the measure of ${T}$.}
\end{defi}

We introduce geometric conditions to obtain the optimal convergence rate of the anisotropic error estimates.

\begin{assume} \label{ass1}
A family of meshes $\{ \mathbb{T}_h\}$ is semi-regular if there exists $\gamma_0 \> 0$ such that
\begin{align}
\displaystyle
\frac{H_{T}}{h_{T}} \leq \gamma_0 \quad \forall \mathbb{T}_h \in \{ \mathbb{T}_h \}, \quad \forall T \in \mathbb{T}_h. \label{NewGeo}
\end{align}
\end{assume}

\begin{rem} \label{goodel}
We consider the good elements on the meshes in \cite{IshKobTsu23a}. On anisotropic meshes, good elements may satisfy the following conditions:
\begin{description}
  \item[($d=2$)] $h_2 \approx h_2 t$;
  \item[($d=3$)] $h_2 \approx h_2 t_1$ and $h_3 \approx h_3 t_2$.
\end{description}
\end{rem}

\begin{rem}
The geometric condition in \eqref{NewGeo}	is equivalent to the maximum angle condition ( \cite[Theorem 1]{IshKobTsu23a}).
\end{rem}

\subsection{Main theorem}
The $L^2$-orthogonal projection $\Pi_{\widehat{T}}^0:L^1(\widehat{T}) \to \mathbb{P}^0(\widehat{T})$ is defined such that, for any $\hat{\varphi} \in L^1(\widehat{T})$,
\begin{align*}
\displaystyle
\Pi_{\widehat{T}}^0 \hat{\varphi} := \frac{1}{|\widehat{T}|} \int_{\widehat{T}}  \hat{\varphi} d\hat{x},
\end{align*}
In setting ${\varphi} := \hat{\varphi} \circ {\Phi}^{-1}$, the associated $L^2$-orthogonal projection $\Pi_{{T}}^0:L^1({T}) \to \mathbb{P}^0({T})$ is defined as
\begin{align*}
\displaystyle
\Pi_{{T}}^0 {\varphi} := \frac{1}{|{T}|} \int_{{T}} {\varphi} d{x},
\end{align*}
where ${T} = {\Phi} (\widehat{T})$. The following theorem provides an anisotropic error estimate for the projection $\Pi_{T}^0$.

\begin{thr} \label{thr1}
{
Let $p \in [1,\infty)$ and $q \in [1,\infty)$ be such that
\begin{align}
\displaystyle
W^{1,p}({T}) \hookrightarrow L^q({T}), \label{Sobolev511}
\end{align}
that is $1 - \frac{d}{p} \geq - \frac{d}{q}$.} Then, for any $\hat{\varphi} \in W^{1,p}(\widehat{T})$ with ${\varphi} := \hat{\varphi} \circ {\Phi}^{-1}$,
\begin{align}
\displaystyle
\| \Pi_{T}^0 \varphi - \varphi \|_{L^q(T)} \leq c |T|^{\frac{1}{q} - \frac{1}{p}} \sum_{i=1}^d h_i \left\| \frac{\partial \varphi}{\partial r_i} \right\|_{L^{p}(T)}. \label{L2ortho}
\end{align}
If Assumption \ref{assume0} is imposed, then:
\begin{align}
\displaystyle
\| \Pi_{T}^0 \varphi - \varphi \|_{L^q(T)} \leq c |T|^{\frac{1}{q} - \frac{1}{p}} \sum_{i=1}^d \widetilde{\mathscr{H}}_i \left\| \frac{\partial (\varphi \circ \Phi_T)}{\partial \tilde{x}_i} \right\|_{L^{p}(\Phi_{T}^{-1}(T))}. \label{L2ortho=ass0}
\end{align}
\end{thr}

\begin{pf*}
The standard scaling argument yields the following:
\begin{align}
\displaystyle
\| \Pi_{T}^0 \varphi - \varphi \|_{L^q(T)}
\leq c |\det({A})|^{\frac{1}{q}} \| \Pi_{\widehat{T}}^0 \hat{\varphi} - \hat{\varphi} \|_{L^q(\widehat{T})}. \label{th1=2}
\end{align}
For any $\hat{\eta} \in \mathbb{P}^{0}$, we have that 
\begin{align}
\displaystyle
\| \Pi_{\widehat{T}}^0 \hat{\varphi} - \hat{\varphi} \|_{L^q(\widehat{T})}
&\leq \| \Pi_{\widehat{T}}^0  (\hat{\varphi} - \hat{\eta}) \|_{L^q(\widehat{T})}  + \| \hat{\eta} - \hat{\varphi} \|_{L^q(\widehat{T})},  \label{th1=3}
\end{align}
because $\Pi_{\widehat{T}}^0 \hat{\eta} = \hat{\eta}$. Using H\"older's inequality yields
\begin{align}
\displaystyle
\| \Pi_{\widehat{T}}^0  (\hat{\varphi} - \hat{\eta}) \|_{L^q(\widehat{T})} 
&= |\widehat{T}|^{\frac{1}{q}} |\Pi_{\widehat{T}}^0 (\hat{\varphi} - \hat{\eta})| 
\leq |\widehat{T}|^{\frac{1}{q} -1} \| \hat{\varphi} - \hat{\eta} \|_{L^1(\widehat{T})} \notag\\
&\leq  |\widehat{T}|^{\frac{1}{q} -\frac{1}{p}} \| \hat{\varphi} - \hat{\eta} \|_{L^p(\widehat{T})}.  \label{th1=4}
\end{align}
Based on \eqref{th1=2}, \eqref{th1=3}, \eqref{th1=4}, and the Sobolev embedding theorem, we have that
\begin{align}
\displaystyle
\| \Pi_{T}^0 \varphi - \varphi \|_{L^q(T)}
\leq  C_1(\widehat{T}) |\det({A})|^{\frac{1}{q}} \inf_{\hat{\eta} \in \mathbb{P}^{0}} \| \hat{\varphi} - \hat{\eta} \|_{W^{1,p}(\widehat{T})}. \label{th1=5}
\end{align}
From the Bramble--Hilbert lemma (refer to \cite[Lemma 4.3.8]{BreSco08}), a constant $\hat{\eta}_{\beta} \in \mathbb{P}^{0}$ exists such that for any $\hat{\varphi} \in W^{1,p}(\widehat{T})$,
\begin{align}
\displaystyle
| \hat{\varphi} - \hat{\eta}_{\beta} |_{W^{s,p}(\widehat{T})} \leq C^{BH}(\widehat{T}) |\hat{\varphi}|_{W^{1,p}(\widehat{T})}, \quad s=0,1. \label{th1=6}
\end{align}
Using the inequality in \cite[Lemma 6]{IshKobTsu23a} with $m=0$, we can estimate inequality \eqref{th1=6} as
\begin{align}
\displaystyle
|\hat{\varphi}|_{W^{1,p}(\widehat{T})}
\leq c |\det({{A}})|^{-\frac{1}{p}} \sum_{i=1}^d h_i \left\| \frac{\partial \varphi}{\partial r_i} \right\|_{L^{p}(T)}.  \label{th1=7}
\end{align}
Using \eqref{th1=5}, \eqref{th1=6}, \eqref{th1=7}, and $ |\det({{A}})| = \frac{|T|}{|\widehat{T}|}$, we can deduce the target inequality \eqref{L2ortho}.

If Assumption \ref{assume0} is imposed, using the inequality in \cite[Lemma 5]{IshKobTsu23a} with $m=0$, we can estimate inequality \eqref{th1=6} as follows:
\begin{align}
\displaystyle
|\hat{\varphi}|_{W^{1,p}(\widehat{T})}
\leq c |\det({{A}})|^{-\frac{1}{p}} \sum_{i=1}^d \widetilde{\mathscr{H}}_i \left\| \frac{\partial \tilde{\varphi}}{\partial \tilde{x}_i} \right\|_{L^{p}(\widetilde{T})}.  \label{th1=7b}
\end{align}
Using \eqref{th1=5}, \eqref{th1=6}, \eqref{th1=7b}, and $ |\det({{A}})| = \frac{|T|}{|\widehat{T}|}$, we can deduce the target inequality \eqref{L2ortho=ass0} with $\tilde{\varphi} = \varphi \circ \Phi_T$ and $\widetilde{T} = \Phi_{T}^{-1}(T)$.
\qed
\end{pf*}

\begin{note} \label{newnote3}
In inequality \eqref{L2ortho=ass0}, obtaining the estimates in $T$ by specifically determining the matrix $A_{T}$ is possible.
Let $p=q=2$. Recall that
\begin{align*}
\displaystyle
\Phi_T: \widetilde{T} \ni \tilde{x} \mapsto x := {A}_{T} \tilde{x} + b_T \in T.
\end{align*}
The space $\mathcal{C}^{1}(\widetilde{T})$ is dense in the space ${H}^{1}(\widetilde{T})$. For $\tilde{\varphi} \in \mathcal{C}^{2}(\widetilde{T})$ with $\varphi = \tilde{\varphi} \circ \Phi_{T}^{-1}$ and $1 \leq i \leq d$, we have
\begin{align*}
\displaystyle
\left| \frac{\partial \tilde{\varphi}}{\partial \tilde{x}_i} (\tilde{x}) \right|
&= \left| \sum_{i_1^{(1)}=1}^d  [A_{T}]_{i_1^{(1)} i} \frac{\partial \varphi}{\partial x_{i_1^{(1)}}} (x)  \right|.
\end{align*}
Let $d=2$. We define a rotation matrix $A_{T}$ as
\begin{align*}
\displaystyle
A_{T} := 
\begin{pmatrix}
\cos \theta  & - \sin \theta \\
 \sin \theta & \cos \theta
\end{pmatrix},
\end{align*}
where $\theta$ denotes the angle. We then have
\begin{align*}
\displaystyle
\left| \frac{\partial \tilde{\varphi}}{\partial \tilde{x}_1} (\tilde{x}) \right|
&=  \left| \cos \theta  \frac{\partial \varphi}{\partial x_{1}} (x)  +  \sin \theta  \frac{\partial \varphi}{\partial x_{2}} (x) \right|, \\
\left| \frac{\partial \tilde{\varphi}}{\partial \tilde{x}_2} (\tilde{x}) \right|
&=  \left| - \sin \theta  \frac{\partial \varphi}{\partial x_{1}} (x)  +  \cos \theta  \frac{\partial \varphi}{\partial x_{2}} (x) \right|.
\end{align*}
If $|\sin \theta| \leq c \frac{\widetilde{\mathscr{H}}_2}{\widetilde{\mathscr{H}}_1}$ and $\widetilde{\mathscr{H}}_2 \leq c \widetilde{\mathscr{H}}_1$, we can deduce
\begin{align*}
\displaystyle
\left| \frac{\partial \tilde{\varphi}}{\partial \tilde{x}_1} (\tilde{x}) \right|
&\leq \left| \frac{\partial \varphi}{\partial x_{1}} (x) \right| + c \frac{\widetilde{\mathscr{H}}_2}{\widetilde{\mathscr{H}}_1} \left| \frac{\partial \varphi}{\partial x_{2}} (x) \right| , \\
\left| \frac{\partial \tilde{\varphi}}{\partial \tilde{x}_2} (\tilde{x}) \right|
&\leq c \left| \frac{\partial \varphi}{\partial x_{1}} (x) \right| +  \left| \frac{\partial \varphi}{\partial x_{2}} (x) \right|.
\end{align*}
As $|\det (A_{T})| = 1$, it holds that for $i=1,2$,
\begin{align*}
\displaystyle
\widetilde{\mathscr{H}}_i \left \| \frac{\partial \tilde{\varphi}}{\partial \tilde{x}_i} \right \|_{L^2(\widetilde{T})}
&\leq c \sum_{j=1}^2 \widetilde{\mathscr{H}}_j  \left \|  \frac{\partial \varphi}{\partial x_{{j}}}  \right \|_{L^2(T)},
\end{align*}
which leads to
\begin{align*}
\displaystyle
\| \Pi_{T}^0 \varphi - \varphi \|_{L^2(T)} \leq c \sum_{j=1}^2 \widetilde{\mathscr{H}}_j \left \|  \frac{\partial \varphi}{\partial x_{{j}}}  \right \|_{L^2(T)}.
\end{align*}

As a special case, we set $\theta:= \frac{\pi}{2}$. We then have
\begin{align*}
\displaystyle
\left| \frac{\partial \tilde{\varphi}}{\partial \tilde{x}_1} (\tilde{x}) \right|
&= \left| \frac{\partial \varphi}{\partial x_{2}} (x)  \right|, \quad \left| \frac{\partial \tilde{\varphi}}{\partial \tilde{x}_2} (\tilde{x}) \right|
= \left| \frac{\partial \varphi}{\partial x_{1}} (x)  \right|.
\end{align*}
As $|\det (A_{T})| = 1$, it holds that
\begin{align*}
\displaystyle
\left \| \frac{\partial \tilde{\varphi}}{\partial \tilde{x}_i} \right \|_{L^2(\widetilde{T})} \leq  \left \|  \frac{\partial \varphi}{\partial x_{{i+1}}}  \right \|_{L^2(T)},
\end{align*}
where indices $i$ and $i+1$ must be evaluated as in Modulo 2. We then have
\begin{align*}
\displaystyle
\| \Pi_{T}^0 \varphi - \varphi \|_{L^2(T)} \leq c \sum_{i=1}^2 \widetilde{\mathscr{H}}_i \left \|  \frac{\partial \varphi}{\partial x_{{i+1}}}  \right \|_{L^2(T)}.
\end{align*}

When $A_T$ is a mirror imaging matrix, similar estimates may hold. Furthermore, the same argument can be made in the case of $d=3$.
\end{note}

\section{Interpolation error estimates of nonconforming finite element methods}
We introduce the following theorem.

\begin{thr} \label{Thr2a}
Let $\alpha := (\alpha_1,\ldots,\alpha_d){^{\top}} \in \mathbb{N}_0^d$ be a multi-index and $k \in \mathbb{N}$. {Let $p \in [1,\infty)$ and $q \in [1,\infty)$ be such that \eqref{Sobolev511} holds.}
We define an interpolation operator $I_{T}: W^{k , p}(T) \to \mathbb{P}^{k}(T)$ that satisfies:
\begin{align}
\displaystyle
\partial^{\alpha} (I_{{T}} {\varphi}) = \Pi_T^0 (\partial^{\alpha} \varphi) \quad \forall \varphi \in W^{k, p}(T) \quad {\forall \alpha: \ |\alpha| \leq k}. \label{sp=cond} 
\end{align}
Then, for any $\hat{\varphi} \in W^{k+1,p}(\widehat{T})$ with ${\varphi} := \hat{\varphi} \circ {\Phi}^{-1}$ and any $\alpha$ with $|\alpha| \leq k$,
\begin{align}
\displaystyle
| I_{T} \varphi - \varphi |_{W^{k,q}({T})} &\leq  c |T|^{ \frac{1}{q} - \frac{1}{p} }  \sum_{i=1}^d  h_i \left | \frac{\partial\varphi}{\partial r_i} \right |_{W^{k,p}(T)}. \label{CR21=43} 
\end{align}
If Assumption \ref{assume0} is imposed, then:
\begin{align}
\displaystyle
| I_{T} \varphi - \varphi |_{W^{k,q}({T})} &\leq  c |T|^{ \frac{1}{q} - \frac{1}{p} }  \sum_{i=1}^d  \widetilde{\mathscr{H}}_i \left | \frac{\partial (\varphi \circ \Phi_T)}{\partial \tilde{x}_i} \right |_{W^{k,p}(\Phi_T^{-1}(T))}. \label{CR21=ass0} 
\end{align}
\end{thr}

\begin{pf*}
The error estimate of the $L^2$-orthogonal projection \eqref{L2ortho} yields
\begin{align*}
\displaystyle
| I_{T} \varphi - \varphi |_{W^{k,q}({T})}^q
&= \sum_{|\alpha| = k} \left \| \partial^{\alpha}  ( I_{T} \varphi - \varphi)  \right\|^q_{L^q(T)} \\
&=  \sum_{|\alpha| = k} \left \| \Pi_T^0 ( \partial^{\alpha} \varphi ) - \partial^{\alpha} \varphi \right\|^q_{L^q(T)} \\
&\leq c |T|^{\left( \frac{1}{q} - \frac{1}{p} \right) q} \sum_{|\alpha| = k} \sum_{i=1}^d  h_i^q \left \| \partial^{\alpha} \frac{\partial \varphi}{\partial r_i} \right \|_{L^p(T)}^q.
\end{align*}
Using the Jensen-type inequality \eqref{jensen}, we obtain
\begin{align*}
\displaystyle
| I_{T} \varphi - \varphi |_{W^{1,q}({T})}
&\leq c  |T|^{ \frac{1}{q} - \frac{1}{p}} \left( \sum_{i=k}^d  h_i^p \sum_{|\alpha| = 1} \left \|  \partial^{\alpha} \frac{\partial\varphi}{\partial r_i} \right \|_{L^p(T)}^p \right)^{\frac{1}{p}} \\
&\leq c  |T|^{ \frac{1}{q} - \frac{1}{p}} \sum_{i=1}^d  h_i \left | \frac{\partial\varphi}{\partial r_i} \right |_{W^{k,p}(T)},
\end{align*}
which is the target inequality in Eq. \eqref{CR21=43}.

If Assumption \ref{assume0} is imposed, then \eqref{L2ortho=ass0} yields
\begin{align*}
\displaystyle
| I_{T} \varphi - \varphi |_{W^{k,q}({T})}^q
&\leq c |T|^{\left( \frac{1}{q} - \frac{1}{p} \right) q} \sum_{|\alpha| = k} \sum_{i=1}^d  \widetilde{\mathscr{H}}_i^q \left \| \partial^{\alpha} \frac{\partial (\varphi \circ \Phi_T)}{\partial \tilde{x}_i} \right \|_{L^p(\Phi_T^{-1}(T))}^q.
\end{align*}
Using the Jensen-type inequality in \eqref{jensen}, we obtained the target inequality in \eqref{CR21=ass0}.
\qed
\end{pf*}

\begin{note}
The operators that satisfy \eqref{sp=cond} exist; see Theorems \ref{thr=CR} and \ref{thr2}.
\end{note}

\subsection{Crouzeix--Raviart interpolation error estimate}
Let ${F}_i$ and $1 \leq i \leq d+1$ be the $(d-1)$-dimensional sub-simplex of ${T}$ opposite vertex ${p}_i$. The CR finite element on the reference element is defined by the triple $\{ {T}, {P}, {\Sigma\}}$ as follows:
\begin{enumerate}
 \item ${P} := \mathbb{P}^1({T})$;
 \item ${\Sigma}$ is a set $\{ {\chi}_{i} \}_{1 \leq i \leq N^{(d,1)}}$ of $N^{(d,1)}$ linear forms $\{ {\chi}_{i} \}_{1 \leq i \leq N^{(d,1)}}$ with its components such that, for any ${q} \in {P}$,
\begin{align}
\displaystyle
{\chi}_{i}({q}) := \frac{1}{| {F}_i |} \int_{{F}_i} {q} d{s} \quad \forall i \in \{ 1, \ldots,d+1 \}. \label{CR911}
\end{align}
\end{enumerate}
Using the barycentric coordinates $ \{ {\lambda}_i \}_{i=1}^{d+1}: \mathbb{R}^d \to \mathbb{R}$ on the reference element, the nodal basis functions associated with the degrees of freedom using \eqref{CR911} is defined as follows:
\begin{align}
\displaystyle
{\theta}_i({x}) := d \left( \frac{1}{d} - {\lambda}_i({x}) \right) \quad \forall i \in \{ 1, \ldots,d+1 \}. \label{CR912}
\end{align}
Thus, ${\chi}_{i} ({\theta}_j) = \delta_{ij}$ for any $i,j \in \{ 1, \ldots,d+1 \}$, The local operator ${I}_{{T}}^{CR}$ is defined as 
\begin{align}
\displaystyle
I_{{T}}^{CR}:  W^{1,1}({T})  \ni {\varphi}  \mapsto I_{{T}}^{CR} {\varphi} := \sum_{i=1}^{d+1} \left( \frac{1}{| {F}_i|} \int_{{F}_i} {\varphi} d{s} \right) {\theta}_i \in {P}. \label{CR913}
\end{align}

We present anisotropic CR interpolation error estimates (see \cite{ApeNicSch01}). 

\begin{thr} \label{thr=CR}
Let $p \in [1,\infty)$ and $q \in [1,\infty)$ be such that \eqref{Sobolev511} holds. Subsequently, for any ${\varphi} \in W^{2,p}({T})$,we have
\begin{align}
\displaystyle
| I_{T}^{CR} \varphi - \varphi |_{W^{1,q}({T})} &\leq  c |T|^{ \frac{1}{q} - \frac{1}{p} }  \sum_{i=1}^d  h_i \left | \frac{\partial\varphi}{\partial r_i} \right |_{W^{1,p}(T)}. \label{CR21} 
\end{align}
If Assumption \ref{assume0} is imposed, then:
\begin{align}
\displaystyle
| I_{T}^{CR} \varphi - \varphi |_{W^{1,q}({T})} &\leq  c |T|^{ \frac{1}{q} - \frac{1}{p} }  \sum_{i=1}^d  \widetilde{\mathscr{H}}_i \left | \frac{\partial (\varphi \circ \Phi_T)}{\partial \tilde{x}_i} \right |_{W^{1,p}(\Phi_T^{-1}(T))}. \label{th3=ass0} 
\end{align}
\end{thr}

\begin{pf*}
Only CR interpolation satisfies the condition \eqref{sp=cond}.

For ${\varphi} \in W^{2,p}({T})$, Green's formula and the definition of the CR interpolation imply that because $I_{T}^{CR} \varphi  \in \mathbb{P}^1$, 
\begin{align*}
\displaystyle
\frac{\partial }{\partial {x}_j} (I_{T}^{CR} \varphi)
&= \frac{1}{|T|} \int_{T}  \frac{\partial }{\partial {x}_j} (I_{T}^{CR} \varphi) dx = \frac{1}{|T|} \sum_{i=1}^{d+1} n_{{T}}^{(j)} \int_{F_i} I_{T}^{CR} \varphi ds \\
&= \frac{1}{|T|} \sum_{i=1}^{d+1} n_{{T}}^{(j)} \int_{F_i} \varphi ds 
= \frac{1}{|T|} \int_{T}  \frac{\partial \varphi}{\partial {x}_j} dx = \Pi_{T}^0 \left(  \frac{\partial \varphi}{\partial {x}_j} \right)
\end{align*}
for  $j=1,\ldots,d$, where $n_{{T}}$ denotes the outer unit normal vector to $T$ and $ n_{{T}}^{(j)}$ denotes the $j$th component of $n_T$. {Therefore, from Theorem \ref{Thr2a}, the target inequalities \eqref{CR21} and \eqref{th3=ass0} hold.}
\qed
\end{pf*}

\subsection{Morley interpolation error estimate}
The Morley FEM has not been defined uniquely. There are two versions: one defined in \cite{Mol68}, which is the original paper, and the other in \cite{ArnBre85,LasLea75,WanJin06}. {In original  Morley FEM, by normal derivatives on faces, the spans of the nodes are not preserved under push-forward. To overcome this difficulty, the mean value of the first normal derivative is used \cite{ArnBre85,LasLea75,WanJin06}. The original Morley interpolation error estimates are obtained using the modified Morley interpolation error estimates (see \cite{LasLea75}). }In this study, we use the Morley FEM introduced in  \cite{WanJin06}.

Let ${F}_i$, $1 \leq i \leq d+1$ be the $(d-1)$-dimensional subsimplex of ${T}$ without vertices  ${p}_i$ and ${S}_{i,j}$, $1 \leq i \< j \leq d+1$ be the $(d-2)$-dimensional subsimplex of ${T}$ without vertices ${p}_i$ and ${p}_j$. The $d$-dimensional Morley finite element on the reference element is defined by the triple $\{ {T}, {P}, {\Sigma\}}$ as 
\begin{enumerate}
 \item ${P} := \mathbb{P}^2({T})$;
 \item ${\Sigma}$ is a set $\{ {\chi}_{i} \}_{1 \leq i \leq N^{(d,2)}}$ of $N^{(d,2)} $ linear forms $\{ {\chi}^{(1)}_{i,j} \}_{1 \leq i \< j \leq d+1} \cup \{ {\chi}^{(2)}_{i} \}_{1 \leq i \leq d+1}$ with its components such that, for any ${q} \in {P}$,
\begin{subequations} \label{Mor1}
\begin{align}
\displaystyle
{\chi}^{(1)}_{i,j}({q}) &:= \frac{1}{|{S}_{i,j}|} \int_{{S}_{i,j}} {q} d{s}, \quad 1 \leq i \< j \leq d+1,  \label{Mor1=a}\\
{\chi}^{(2)}_{i}({q}) &:= \frac{1}{|{F}_{i}|} \int_{{F}_{i}} \frac{\partial {q}}{\partial {n}_i} d{s}, \quad 1 \leq i \leq d+1, \label{Mor1=b}
\end{align}
\end{subequations}
where $\frac{\partial }{\partial {n}_i} = n_{{T,i}} \cdot \nabla $, and $n_{{T,i}}$ is the unit outer normal to ${F}_i \subset \partial {T}$. For $d=2$, ${\chi}^{(1)}_{i,j}({q})$ is interpreted as
\begin{align*}
\displaystyle
{\chi}^{(1)}_{i,j}({q}) = {q}({p}_k), \quad k=1,2,3, \quad k \neq i,j.
\end{align*} 
\end{enumerate}
For a Morley finite {element}, ${\Sigma}$ is unisolvent (see \cite[Lemma 2]{WanJin06}). The nodal basis functions associated with the degrees of freedom provided by \eqref{Mor1} are defined as follows:
\begin{subequations} \label{Mor2}
\begin{align}
\displaystyle
{\theta}_{i,j}^{(1)} &:= 1 - (d-1)({\lambda}_i + {\lambda}_j) + d(d-1) {\lambda}_i {\lambda}_j \notag \\
&\quad - (d-1) (\nabla {\lambda}_i){^{\top}} \nabla {\lambda}_j \sum_{k = i,j} \frac{{\lambda}_k (d {\lambda}_k - 2)}{2 | \nabla {\lambda}_k |_E^2}, \quad 1 \leq i \< j \leq d+1, \label{Mor2=a} \\
{\theta}_{i}^{(2)} &:= \frac{{\lambda}_i (d {\lambda}_i - 2)}{2 | \nabla {\lambda}_i |_E}, \quad 1 \leq i \leq d+1, \label{Mor2=b}
\end{align}
\end{subequations}
where $|\nabla {\lambda}_i|_E$ denotes the Euclidean norm in $\mathbb{R}^d$. Subsequently, \cite[Theorem 1]{WanJin06} proved that, for $1 \leq i \< j \leq d+1$,
\begin{align}
\displaystyle
{\chi}^{(1)}_{k,\ell}({\theta}_{i,j}^{(1)}) = \delta_{ik} \delta_{j \ell}, \quad 1 \leq k \< \ell \leq d+1,  \quad {\chi}^{(2)}_{k}({\theta}_{i,j}^{(1)}) = 0, \quad 1 \leq k \leq d+1, \label{Mor3}
\end{align}
and, for $1 \leq i \leq d+1$,
\begin{align}
\displaystyle
{\chi}^{(1)}_{k,\ell}({\theta}^{(2)}_{i}) = 0, \quad 1 \leq k \< \ell \leq d+1, \quad {\chi}^{(2)}_{k}({\theta}^{(2)}_{i}) = \delta_{ik}, \quad 1 \leq k \leq d+1. \label{Mor4}
\end{align}
The local interpolation operator ${I}_{{T}}^{M}$ is defined as 
\begin{align}
\displaystyle
{I}_{{T}}^{M}: W^{2,1}({T}) \ni {\varphi} \mapsto {I}_{{T}}^{M} {\varphi} \in {P}, \label{Mor5}
\end{align}
with
\begin{align}
\displaystyle
{I}_{{T}}^{M} {\varphi} := \sum_{1 \leq i \< j \leq d+1} {\chi}^{(1)}_{i,j}({\varphi}) {\theta}_{i,j}^{(1)} + \sum_{1 \leq i \leq d+1} {\chi}^{(2)}_{i}({\varphi}) {\theta}_{i}^{(2)}.  \label{Mor6}
\end{align}
Then, it holds that ${I}_{{T}}^{M} {q} = {q}$ for any ${q} \in {P}$ and ${\varphi} \in W^{2,1}({T})$.
\begin{subequations} \label{Mor7}
\begin{align}
\displaystyle
{\chi}^{(1)}_{i,j}({I}_{{T}}^{M} {\varphi}) &= {\chi}^{(1)}_{i,j}( {\varphi}), \quad 1 \leq i \< j \leq d+1, \label{Mor7=a}\\
{\chi}^{(2)}_{i}({I}_{{T}}^{M} {\varphi}) &= {\chi}^{(2)}_{i}( {\varphi}), \quad 1 \leq i \leq d+1.  \label{Mor7=b}
\end{align}
\end{subequations}

{Using the idea of \cite[Lemma 1]{WanJin06}, the following lemma holds.}


\begin{lem} \label{lem=new1}
Let $T \subset \mathbb{R}^d$ be a simplex. $n_{T,k}$ denotes the unit outer normal to the face $F_k$, $k  =1,\ldots,d+1$ of $T$, $S_1,\ldots,S_d$ are all $(d-2)$-dimensional subsimplexes of $F_k$. Let $v \in \mathcal{C}^1(T)$ be such that
\begin{align}
\displaystyle
\int_{S_{\ell}} v = 0, \quad \int_{F_k} \frac{\partial v}{\partial n_k} = 0, \label{Mor22}
\end{align}
for any $\ell=1,\ldots,d$ and  $k=1,\ldots,d+1$. It then holds that
\begin{align}
\displaystyle
 \int_{F_k} \frac{\partial v}{\partial x_i} = 0, \quad i=1,\ldots,d, \quad k=1,\ldots,d+1.  \label{Mor23}
\end{align}
\end{lem}

\begin{pf*}
 Let $v \in \mathcal{C}^1(T)$. Let $\xi \in \mathbb{R}^d$ be a constant vector, and let $\tau := \xi - (\xi \cdot n_{T,k}) n_{T,k}$. We have
\begin{align*}
\displaystyle
\tau \cdot n_{T,k}
&= \xi \cdot n_{T,k} - (\xi \cdot n_{T,k}) n_{T,k} \cdot n_{T,k} = 0,
\end{align*}
that is, $\tau$ is the tangent vector of $F_k$. Subsequently, from \eqref{Mor22} we obtain
\begin{align*}
\displaystyle
\int_{F_k} (\xi \cdot \nabla ) v 
= \int_{F_k}  \frac{\partial v}{\partial \tau} +  (\xi \cdot n_{T,k}) \int_{F_k} \frac{\partial v}{\partial n_k}
=  \int_{F_k}  \frac{\partial v}{\partial \tau}.
\end{align*}
Let $d=2$. Let $p_{k1}$ and $p_{k2}$ be the endpoints of the edge $F_k$, that is, $F_k = \overline{p_{k1} p_{k2}}$. Subsequently,  from \eqref{Mor22} we obtain
\begin{align}
\displaystyle
 \int_{F_k}  \frac{\partial v}{\partial \tau}
 &= \int_{s=0}^{s=|F_k|} \frac{d v}{ds} \left( \frac{|F_k| - s}{|F_k|}p_{k1} + \frac{s}{|F_k|} p_{k2} \right) = v(p_{k2}) - v(p_{k1}) = 0. \label{Mor21=2d}
\end{align}
Let $d=3$.  Let $\zeta^{(\ell)}$ be the unit outer normal of $S_{\ell}$ for $\ell=1,2,3$, From \eqref{Mor22}, the Gauss--Green formula yields
\begin{align}
\displaystyle
 \int_{F_k}  \frac{\partial v}{\partial \tau}
&= \sum_{\ell=1}^3 \tau \cdot \zeta^{(\ell)} \int_{S_{\ell}} v = 0. \label{Mor21}
\end{align}
From \eqref{Mor21=2d} and \eqref{Mor21}, it holds that for $d=2,3$
\begin{align}
\displaystyle
\int_{F_k} (\xi \cdot \nabla ) v = 0. \label{Mor21=23d}
\end{align}
Let $e_1, \ldots, e_d \in \mathbb{R}^d$ be a canonical basis. By setting $\xi := e_i$ in \eqref{Mor21=23d}, we obtain the desired result in \eqref{Mor23} under Assumption \eqref{Mor22}.
\qed
\end{pf*}

The anisotropic Morley interpolation error estimate is expressed as 

\begin{thr} \label{thr2}
Let $p \in [1,\infty)$ and $q \in [1,\infty)$ be such that \eqref{Sobolev511} holds true. Subsequently, for any ${\varphi} \in W^{3,p}({T}) \cap \mathcal{C}^1({T})$,we have
\begin{align}
\displaystyle
| I_{T}^{M} \varphi - \varphi |_{W^{2,q}({T})} &\leq  c |T|^{ \frac{1}{q} - \frac{1}{p} }  \sum_{i=1}^d  h_i \left | \frac{\partial\varphi}{\partial r_i} \right |_{W^{2,p}(T)}. \label{Mor24} 
\end{align}
If Assumption \ref{assume0} is imposed, then:
\begin{align}
\displaystyle
| I_{T}^{M} \varphi - \varphi |_{W^{2,q}({T})} &\leq  c |T|^{ \frac{1}{q} - \frac{1}{p} }  \sum_{i=1}^d  \widetilde{\mathscr{H}}_i \left | \frac{\partial (\varphi \circ \Phi_T)}{\partial \tilde{x}_i} \right |_{W^{2,p}(\Phi_T^{-1}(T) )}. \label{Mor24=ass0} 
\end{align}
\end{thr}

\begin{pf*}
Only Morley interpolation satisfies the condition \eqref{sp=cond}.

Let ${\varphi} \in W^{3,p}({T})  \cap \mathcal{C}^1(T)$ and set $v :=  I_{T}^{M} \varphi - \varphi$. Using the definition of the Morley interpolation operator \eqref{Mor7}, we obtain
\begin{align*}
\displaystyle
 \int_{{S}_{i,j}} {v} d {s} = 0, \quad 1 \leq i \< j \leq d+1, \quad
 \int_{{F}_{i}} \frac{\partial v}{ \partial n_i} ds = 0, \quad 1 \leq i \leq d+1.
\end{align*}
Therefore, from Lemma \ref{lem=new1}, we have
\begin{align}
\displaystyle
\int_{F_i} \frac{\partial v}{\partial x_k} = 0, \quad i=1,\ldots,d+1, \quad k = 1,\ldots,d. \label{Mor26}
\end{align}
From Green's formula and \eqref{Mor26}, it follows that, for $1 \leq j,k \leq d$,
\begin{align*}
\displaystyle
\int_{{T}} \frac{\partial^2 v}{\partial {x}_j \partial {x}_k} d {x}
&= \sum_{i=1}^{d+1}  n_{{T},j} \int_{{F}_i}  \frac{\partial v}{\partial {x}_k} d {s} = 0,
\end{align*}
which leads to
\begin{align*}
\displaystyle
 \frac{\partial^2}{\partial {x}_j \partial {x}_k} (I_{{T}}^{M} {\varphi} )
  = \frac{1}{|T|} \int_{T}  \frac{\partial^2 }{\partial {x}_j \partial {x}_k}  (I_{{T}}^{M} {\varphi} ) dx
 = \frac{1}{|T|} \int_{T}  \frac{\partial^2 \varphi}{\partial {x}_j \partial {x}_k} dx = \Pi_{T}^0 \left( \frac{\partial^2 \varphi}{\partial {x}_j \partial {x}_k} \right),
\end{align*}
because $ \frac{\partial^2}{\partial {x}_j \partial {x}_k} (I_{{T}}^{M} {\varphi}) \in \mathcal{P}^0(T)$, {Therefore, by Theorem \ref{Thr2a}, the target inequalities \eqref{Mor24} and \eqref{Mor24=ass0} hold.}
\qed
\end{pf*}

\section{Preliminaries for anisotropic error analysis}

\subsection{Useful relation for anisotropic error analysis}
For any $T \in \mathbb{T}_h$, the local RT polynomial space is defined as 
\begin{align*}
\displaystyle
\mathbb{RT}^0(T) := \mathbb{P}^0(T)^d + x \mathbb{P}^0(T), \quad x \in \mathbb{R}^d.
\end{align*}
The RT finite element space is defined as follows:
\begin{align*}
\displaystyle
V^{RT^0}_{h} &:= \{ v_h \in L^1(\Omega)^d: \  v_h|_T \in {\mathbb{RT}^0(T)}, \ \forall T \in \mathbb{T}_h, \  [\![ v_h \cdot n ]\!]_F = 0, \ \forall F \in \mathcal{F}_h \}.
\end{align*}
The discontinuous and CR finite element spaces are defined as
\begin{align*}
\displaystyle
P_{dc,h}^m &:= \left\{ p_h \in L^2(\Omega); \ p_h|_{T} \in \mathbb{P}^{m}({T}) \quad \forall T \in \mathbb{T}_h \right\}, \quad m \in \mathbb{N}_0, \\
V_{h0}^{CR} &:=  \biggl \{ \varphi_h \in P_{dc,h}^1: \  \int_F [\![ \varphi_h ]\!] ds = 0 \ \forall F \in \mathcal{F}_h \biggr \}
\end{align*}
and the Morley finite element space is as follows:
\begin{align*}
\displaystyle
V_{h}^M &:=  \biggl \{ \varphi_h \in P_{dc,h}^2: \ \int_F  \biggl {[} \! \biggl {[} \frac{ \partial \varphi_h}{  \partial n}  \biggr {]} \! \biggr {]} ds = 0 \ \forall F \in \mathcal{F}_h^i, \\
&\hspace{0.7cm}  \text{the integral average  of $\varphi_h$ over each $(d-2)$-dimensional} \\
&\hspace{0.7cm} \text{subsimplex of $T \in \mathbb{T}_h$ is continuous} \biggr \}, \\
V_{h0}^M &:=  \left\{ \varphi_h \in V_{h}^M; \ \text{degrees of freedom of $\varphi_h$ in \eqref{Mor1} vanish on $\partial \Omega$} \right \}.
\end{align*}
In particular, for $d=2$, the space $V_{h0}^M$ is described as
\begin{align*}
\displaystyle
V_{h0}^M &:=  \biggl \{ \varphi_h \in P_{dc,h}^2: \ \int_F  \biggl {[} \! \biggl {[} \frac{ \partial \varphi_h}{  \partial n}  \biggr {]} \! \biggr {]} ds = 0 \ \forall F \in \mathcal{F}_h, \\
&\hspace{0.7cm} \text{$\varphi_h$ is continuous at each vertex in $\Omega$, $\varphi_h(p) = 0$, \ $p \in \partial \Omega$} \biggr \}.
\end{align*}

The following relationship is crucial in the error analysis of nonconforming finite elements on anisotropic meshes and holds for isotropic meshes:

\begin{lem} \label{lem3}
For any $v_h \in V^{RT^0}_{h}$ or $\psi_h \in H_0^1(\Omega) +  V_{h0}^{CR}$,we have
\begin{align}
\displaystyle
\sum_{T \in \mathbb{T}_h} \int_{T} ( v_h \cdot \nabla ) \psi_h dx + \sum_{T \in \mathbb{T}_h} \int_{T} \div v_h \psi_h dx  = 0,  \label{CR61}
\end{align}
For any $v_h \in V^{RT^0}_{h}$ or $\varphi_h \in V_{h0}^M$
\begin{align}
\displaystyle
\sum_{T \in \mathbb{T}_h} \int_{T} \left( v_h \cdot \nabla \right) \frac{\partial \varphi_h}{\partial x_i} dx + \sum_{T \in \mathbb{T}_h}  \int_{T} \div v_h  \frac{\partial \varphi_h}{\partial x_i} dx = 0, \quad i = 1,\ldots,d.  \label{Mh=5}
\end{align}
\end{lem}

\begin{pf*}
For any $v_h \in V_h^{RT^0}$ and $\psi_h \in H_0^1(\Omega) +  V_{h0}^{CR}$, using Green’s formula and that $v_h \cdot n_F \in \mathbb{P}^0(F)$ for any $F \in \mathcal{F}_h$, we can derive 
\begin{align*}
\displaystyle
&\sum_{T \in \mathbb{T}_h} \int_{T} ( v_h \cdot \nabla ) \psi_h dx + \sum_{T \in \mathbb{T}_h} \int_{T} \div v_h \psi_h dx \\
&\quad = \sum_{T \in \mathbb{T}_h} \int_{\partial T} (v_h \cdot n_{T}) \psi_h ds  \notag 
 = \sum_{F \in \mathcal{F}_h} \int_{F} [\![ (v_h  \psi_h ) \cdot n_F ]\!] ds  \notag \\
&\quad = \sum_{F \in \mathcal{F}_h^i} \int_{F} \left(  [\![ v_h \cdot n_F]\!] \{ \! \{ \psi_h \} \!\} + \{ \! \{ v_h  \} \!\} \cdot n_F [\![ \psi_h ]\!] \right) ds \notag \\
&\quad \quad + \sum_{F \in \mathcal{F}_h^{\partial}} \int_{F} (v_h \cdot n_F) \psi_h ds  =  0,
\end{align*}
and are given by \eqref{CR61}.

Let $\varphi_h \in V_{h0}^M$. Suppose that $F \in \mathcal{F}_h^i$ with $F = T_{\natural} \cap T_{\sharp}$, $T_{\natural},T_{\sharp} \in \mathbb{T}_h$. Let $S_1,\ldots,S_d$ denote all $(d-2)$-dimensional subsimplexes of $F$. According to Lemma \ref{lem=new1}, because
\begin{align*}
\displaystyle
\int_{S_{\ell}} [\![ \varphi_h ]\!] ds = 0, \quad \ell \in \{1,\ldots,d \}, \quad \int_F \biggl {[} \! \biggl {[} \frac{ \partial \varphi_h}{  \partial n}  \biggr {]} \! \biggr {]}_{F} ds = 0
\end{align*}
it holds that
\begin{align}
\displaystyle
\int_{F}  \biggl {[} \! \biggl {[} \frac{ \partial \varphi_h}{  \partial x_i}  \biggr {]} \! \biggr {]}_{F} ds = 0 \quad \forall F \in \mathcal{F}_h^i, \ i=1,\ldots,d. \label{Mh=6}
\end{align}
Similarly, it holds that
\begin{align}
\displaystyle
\int_{F}  \frac{\partial \varphi_h}{  \partial x_i}  ds = 0 \quad \forall F \in \mathcal{F}_h^{\partial}, \ i=1,\ldots,d. \label{Mh=7}
\end{align}
By using Green’s formulas \eqref{Mh=6}, \eqref{Mh=7}, and $v_h \cdot n_F \in \mathbb{P}^{0}(F)$ for any $F \in \mathcal{F}_h$, we can derive 
\begin{align*}
\displaystyle
&\sum_{T \in \mathbb{T}_h} \int_{T} ( v_h \cdot \nabla ) \frac{\partial \varphi_h}{\partial x_i} dx + \sum_{T \in \mathbb{T}_h}  \int_{T} \div v_h  \frac{\partial \varphi_h}{\partial x_i} dx \\
&\quad = \sum_{T \in \mathbb{T}_h} \int_{\partial T} (v_h \cdot n_T)  \frac{\partial \varphi_h}{\partial x_i} ds \\
&\quad =  \sum_{F \in \mathcal{F}_h^i} \int_{F} \left(  [\![ v_h \cdot n_F]\!]  \biggl {\{} \! \biggl {\{} \frac{ \partial \varphi_h}{  \partial x_i}  \biggr {\}} \! \biggr {\}}_{F} + \{\! \{ v_h  \}\! \} \cdot n_F  \biggl {[} \! \biggl {[} \frac{ \partial \varphi_h}{  \partial x_i}  \biggr {]} \! \biggr {]}_{F}  \right) ds \\
&\quad \quad + \sum_{F \in \mathcal{F}_h^{\partial}}  \int_{F} ( v_h \cdot n_F ) \frac{ \partial \varphi_h}{  \partial x_i} ds = 0,
\end{align*}
and are given by \eqref{Mh=5}:
\qed	
\end{pf*}

\subsection{Error estimate of the RT finite element method} \label{sec=RT}
Error estimates of RT finite element interpolation on anisotropic meshes can be found in \cite{Acoetal11,AcoDur99,DurLom08,FarNicPaq01,Ish22b}. Here, we introduce our version \cite{Ish22b} using a new geometric parameter (Definition \ref{defi1}).

For $q \in \mathbb{RT}^0(T)$, the local degrees of freedom are defined as 
\begin{align*}
\displaystyle
{\chi}^{RT^0}_{T,i}({q}) := \int_{{F}_{i}} {q} \cdot n_{T,i} d{s}, \quad F_i \subset \partial T, \quad \forall i \in \{ 1, \ldots , d+1 \}.
\end{align*}
When setting $\Sigma_{T}^{RT^0} := \{ {\chi}^{RT^0}_{T,i} \}_{1 \leq i \leq d+1}$, the triple $\{ T,  \mathbb{RT}^0(T) , \Sigma_{T}^{RT^0} \}$ is a finite element. The local shape functions are as follows. 
\begin{align*}
\displaystyle
\theta_{T,i}^{RT^0}(x) := \frac{\iota_{{F_{i}},{T}}}{d |T|} (x - p_{i}) \quad \forall i \in \{ 1, \ldots , d+1 \},
\end{align*}
where $\iota_{{F_{i}},{T}} := 1$ if ${n}_{T,{i}}$ points outwards, and $ - 1$ otherwise \cite[Chapter 14]{ErnGue21a}. 

Let $\mathcal{I}_{T}^{RT^0}: H^{1}(T)^d \to \mathbb{RT}^0(T)$ be the RT interpolation operator such that for any $v \in H^{1}(T)^d$,
\begin{align*}
\displaystyle
\mathcal{I}_{T}^{RT^0}: H^{1}(T)^d \ni v \mapsto \mathcal{I}_{T}^{RT} v := \sum_{i=1}^{d+1} \left(  \int_{{F}_{i}} {v} \cdot n_{T,i} d{s} \right) \theta_{T,i}^{RT^0} \in \mathbb{RT}^0(T).
\end{align*}

The Piola transformation $\Psi : L^1(\widehat{T})^d \to L^1({T})^d$ is defined as follows:
\begin{align*}
\displaystyle
\Psi :  L^1(\widehat{T})^d  &\to  L^1({T})^d \\
\hat{v} &\mapsto v(x) :=  \Psi_{T} \circ {\Psi}_{\widetilde{T}}(\hat{v})(x) = \frac{1}{\det(A)} A \hat{v}(\hat{x}), \quad A = {A}_{T} {A}_{\widetilde{T}}
\end{align*}
with	 the following two Piola transformations:
\begin{align*}
\displaystyle
{\Psi}_{\widetilde{T}}:  L^1(\widehat{T})^d &\to L^1(\widetilde{T})^d \nonumber\\
\hat{v} &\mapsto \tilde{v}(\tilde{x}) := {\Psi}_{\widetilde{T}}(\hat{v})(\tilde{x}) := \frac{1}{\det ({A}_{\widetilde{T}})} {A}_{\widetilde{T}}  \hat{v}(\hat{x}), \\
{\Psi}_{T}:  L^1(\widetilde{T})^d &\to  L^1({T})^d \nonumber\\
\tilde{v} &\mapsto {v}({x}) := {\Psi}_{T}(\tilde{v})({x}) := \frac{1}{\det ({A}_{T})} {A}_{T} \tilde{v}(\tilde{x}).
\end{align*}

The following two theorems are divided into the elements of (Type \roman{sone}) and (Type \roman{stwo}) in Section \ref{sec322} for $d=3$. 


\begin{thr} \label{thr3}
Let $T$ be an element with conditions \ref{cond1} and \ref{cond2} which satisfy (Type \roman{sone}) in Section \ref{sec322} when $d=3$. For any $\hat{v} \in H^{1}(\widehat{T})^d$, ${v} = ({v}_1,\ldots,{v}_d){^{\top}} := {\Psi} \hat{v}$
\begin{align}
\displaystyle
\| \mathcal{I}_{T}^{RT^0} v - v \|_{L^2(T)^d} 
&\leq  c \left( \frac{H_{T}}{h_{T}} \sum_{i=1}^d h_i \left \|  \frac{\partial v}{\partial r_i} \right \|_{L^2(T)^d} +  h_{T} \| \div {v} \|_{L^{2}({T})} \right). \label{RT5}
\end{align}
If Assumption \ref{assume0} is imposed, then:
\begin{align}
\displaystyle
\| \mathcal{I}_{T}^{RT^0} v - v \|_{L^2(T)^d} 
&\leq  c \Biggl( \frac{H_{T}}{h_{T}} \sum_{i=1}^d \widetilde{\mathscr{H}}_i \left \|  \frac{\partial (\Psi_T^{-1} v)}{\partial \tilde{x}_i} \right \|_{L^2(\Phi_T^{-1}(T))^d} \nonumber\\
&\hspace{2cm} +  h_{T} \| \div {(\Psi_T^{-1} v)} \|_{L^{2}(\Phi_T^{-1}(T))} \Biggr). \label{RT5=ass0}
\end{align}
\end{thr}

\begin{pf*}
The proof can be found in \cite[Theorem 2]{Ish22b}.
\qed
\end{pf*}

\begin{thr} \label{newthr6}
Let $d=3$. Let $T$ be an element with Condition \ref{cond2} that satisfies (Type \roman{stwo}) in Section \ref{sec322}. For any $\hat{v} \in H^{1}(\widehat{T})^3$, ${v} = ({v}_1,v_2,{v}_3){^{\top}} := {\Psi} \hat{v}$
\begin{align}
\displaystyle
\| \mathcal{I}_{T}^{RT^0} v - v \|_{L^2(T)^3} 
&\leq  c  \frac{H_{T}}{h_{T}} \left( h_T \sum_{i=1}^3 \left \|  \frac{\partial v}{\partial r_i} \right \|_{L^2(T)^3} \right). \label{newRT5}
\end{align}
If Assumption \ref{assume0} is imposed, then:
\begin{align}
\displaystyle
\| \mathcal{I}_{T}^{RT^0} v - v \|_{L^2(T)^3} 
&\leq  c \frac{H_{T}}{h_{T}} \Biggl( \sum_{i=1}^3 \widetilde{\mathscr{H}}_i \left \|  \frac{\partial (\Psi_T^{-1} v)}{\partial \tilde{x}_i} \right \|_{L^2(\Phi_T^{-1}(T))^3} \nonumber \\
&\hspace{2cm} + h_T \sum_{i=1}^3 \left \|  \frac{\partial (\Psi_T^{-1} v)}{\partial \tilde{r}_i} \right \|_{L^2(\Phi_T^{-1}(T))^3} \Biggr). \label{newRT5=ass0}
\end{align}
\end{thr}

\begin{pf*}
The proof can be found in \cite[Theorem 3]{Ish22b}.
\qed
\end{pf*}

\begin{note} \label{newnote5}
Let $d=3$. An anisotropic RT interpolation error estimate cannot be obtained for (Type \roman{stwo}) in Section \ref{sec322}. Therefore, we {do} not obtain the advantage of using anisotropic meshes. 
\end{note}

We define the global RT interpolation $\mathcal{I}_{h}^{RT^0} : (V_{h0}^{CR})^d  \cup H^1(\Omega)^d \to V^{RT^0}_{h}$ as follows:
\begin{align*}
\displaystyle
(\mathcal{I}_{h}^{RT^0} v )|_{T} = \mathcal{I}_{T}^{RT^0} (v|_{T}) \quad \forall T \in \mathbb{T}_h, \quad \forall v \in  (V_{h0}^{CR})^d \cup H^1(\Omega)^d.
\end{align*}
Furthermore, we define the global interpolation $\Pi_h^0$ to space $P_{dc,h}^{0}$ as:
\begin{align*}
\displaystyle
(\Pi_h^0 \varphi)|_{T} := \Pi_{T}^0 (\varphi|_{T}) \quad \forall T \in \mathbb{T}_h, \quad \forall \varphi \in L^2(\Omega).
\end{align*}
The following relationship holds between the RT interpolation $\mathcal{I}_{h}^{RT^0}$ and $L^2$-projection $\Pi_h^0$:

\begin{lem} \label{lem2}
It holds that
\begin{align}
\displaystyle
\div (\mathcal{I}_{T}^{RT^0} v) = \Pi_{T}^0 (\div v) \quad \forall v \in H^{1}(T_j)^d. \label{RT53}
\end{align}
Incorporating \eqref{RT53}, for any $v \in H^1(\Omega)^d$, it holds that
\begin{align}
\displaystyle
\div (\mathcal{I}_{h}^{RT^0} v) = \Pi_h^0 (\div v). \label{RT54}
\end{align}
\end{lem}

\begin{pf*}
The proof is provided in \cite[Lemma 16.2]{ErnGue21a}.
\qed
\end{pf*}

\subsection{Error estimates of the Lagrange finite element interpolation for $d=2$} \label{sec5=5}
We introduce anisotropic Lagrange finite-element interpolation error estimates for the analysis. (See \cite{IshKobTsu23a} for detailed results.)

Let ${I}_{T}^{L}: \mathcal{C}(T) \to \mathbb{P}^1(T)$ be the Lagrange interpolation operator such that for any $\varphi \in \mathcal{C}(T) $,
\begin{align*}
\displaystyle
{I}_{T}^{L}: \mathcal{C}(T)  \ni \varphi \mapsto I_{T}^{L} \varphi := \sum_{i=1}^{3} \varphi(p_{i}) \theta_{T,i}^{L} \in \mathbb{P}^1(T),
\end{align*}
where the local shape functions are 
\begin{align*}
\displaystyle
\theta_{T,i}^{L}(x) :=  {\lambda}_{T,i}(x), \quad i \in \{ 1 ,2,3 \}.
\end{align*}

\begin{thr} \label{thr4}
For all $\hat{\varphi} \in H^{2}(\widehat{T})$ with ${\varphi} := \hat{\varphi} \circ {\Phi}^{-1}$, we have
\begin{align}
\displaystyle
| {\varphi} - I_{{T}}^L {\varphi}|_{H^{m}({T})}
\leq  c \left( \frac{H_{T}}{h_{T}} \right)^m \sum_{|\varepsilon| =  2-m} {h}^{\varepsilon} | \partial_r^{\varepsilon} \varphi |_{H^{m}(T)}, \quad m=0,1. \label{Lag110}
\end{align}

\end{thr}

\begin{pf*}
The proof can be found in \cite[Corollary 1]{IshKobTsu23a}.
\qed
\end{pf*}

The conforming finite element spaces are defined as
\begin{align*}
\displaystyle
V^{L}_{h} &:= \{ v_h \in H^1(\Omega); v_h |_T \in \mathbb{P}^1(T), \ \forall T \in \mathbb{T}_h \}, \quad V^{L}_{h0} := V^{L}_{h} \cap H_0^1(\Omega).
\end{align*}
We set $N := \dim V^{L}_{h}$. The global interpolation operator is defined as follows: 
\begin{align*}
\displaystyle
I_h^{ML}: V_{h0}^{M} \cup H^2_0(\Omega) \ni \varphi_h \mapsto I_h^{ML} \varphi_h := \sum_{i=1}^N \varphi_h(P_i) \theta_i \in V^{L}_{h0}.
\end{align*}
Let $\{ \mathbb{T}_h\}$ be a family of conformal meshes with semiregular properties (Assumption \ref{ass1}). Subsequently, for all $\varphi_h \in V_{h0}^{M}$
\begin{align}
\displaystyle
| I_h^{ML} \varphi_h - \varphi_h |_{H^1(\mathbb{T}_h)} &\leq c h |\varphi_h|_{H^2(\mathbb{T}_h)},\label{error=22} \\
\| I_h^{ML} \varphi_h - \varphi_h \|_{L^2(\Omega)} &\leq c h^2 |\varphi_h|_{H^2(\mathbb{T}_h)}.\label{error=23}
\end{align}

\section{Application to the fourth-order elliptic problem}
Below, we use the interpolation error estimates in the $r_i$-coordinate system for the analysis.

\subsection{Continuous problem}
Let $\Omega \subset \mathbb{R}^2$ be {a} bounded polygonal domain. The fourth-order elliptic problem involves determining $\psi: \Omega \to \mathbb{R}$ such that:
\begin{align}
\displaystyle
\nu \varDelta^2 \psi &= g   \quad \text{in $\Omega$}, \quad \psi = \frac{\partial \psi}{\partial n} = 0 \quad \text{on $\partial \Omega$} \label{cont=1},
\end{align}
where $g:\Omega \to \mathbb{R}$ is a given function and $\nu$ is a real parameter with $0 \< \nu \leq 1$. We define the bilinear form $a_1: H^{2}(\Omega) \times H^2(\Omega) \to \mathbb{R}$ as follows:
\begin{align*}
\displaystyle
a_1(\psi , \varphi) := \sum_{i=1}^2 \int_{\Omega} \nabla \frac{\partial \psi}{\partial x_i} \cdot \nabla \frac{\partial \varphi}{\partial x_i} dx
\end{align*}
for any $\psi,\varphi \in H^2(\Omega)$. The typical variational formulation of the fourth-order {elliptic} problem \eqref{cont=1} is as follows. For any $g \in H^{-1}(\Omega)$, we determine $\psi \in H_0^2(\Omega)$ such that
\begin{align}
\displaystyle
\nu a_1(\psi,\varphi) = \int_{\Omega} g  \varphi dx \quad \forall \varphi \in H_0^2(\Omega). \label{cont=2} 
\end{align}
{According to the Lax--Milgram theorem, a unique solution $\psi \in H^2_0(\Omega)$ of the problem \eqref{cont=2} exists.} The proof is straightforward (see \cite[pp. 301-302]{Gri11}).

The following theorem is known.
\begin{thr}
For any $g\in H^{-1}(\Omega)$, there exists $\varepsilon \> 0$ such that $u \in H^{\frac{5}{2} + \varepsilon}(\Omega)$. Furthermore, if $\Omega$ is convex, then solution $u$ to problem \eqref{cont=2} belongs to $H^3(\Omega)$. In other words, $\varDelta^2$ is an isomorphism between $H^3(\Omega) \cap H^2_0(\Omega)$ and $H^{-1}(\Omega)$. 
\end{thr}

\begin{pf*}
 The proofs are provided in \cite[Corollaries 3.4.2, 3.4.3]{Gri92} and \cite[Corollary 7.3.2.5]{Gri11}.
\qed
\end{pf*}

\subsection{Modified Morley method for the fourth-order elliptic problem}
For any $g \in L^2(\Omega)$, the typical Morley FEM for problem \eqref{cont=2} involves determining $\psi_h \in V_{h0}^M$ such that
\begin{align}
\displaystyle
\nu a_{h1}(\psi_h , \varphi_h) =  \int_{\Omega} g \varphi_h dx \quad \forall \varphi_h \in V_{h0}^M, \label{mod=1}
\end{align}
where the bilinear form $a_{h1}:( H^{2}(\Omega) + V_{h0}^M) \times (H^2(\Omega) + V_{h0}^M) \to \mathbb{R}$ is:
\begin{align*}
\displaystyle
 a_{h1}(\psi_h , \varphi_h) := \sum_{i=1}^2 \sum_{T \in \mathbb{T}_h} \int_{T} \nabla \frac{\partial \psi_h}{\partial x_i} \cdot \nabla \frac{\partial \varphi_h}{\partial x_i} dx.
\end{align*}
However, because $V_{h0}^M \not\subset H_0^1(\Omega)$, the typical Morley {finite element problem \eqref{mod=1}} cannot be applied to $g \in H^{-1}(\Omega)$. Hence, we consider the modified Morley finite element approximation problem of \eqref{cont=2} proposed in \cite{ArnBre85}. The problem involves determining $\psi_h^{*} \in V_{h0}^M$ such that
\begin{align}
\displaystyle
\nu a_{h1}(\psi_h^{*},\varphi_h) = \int_{\Omega} g  I_h^{ML}  \varphi_h dx \quad \forall \varphi_h \in V^M_{h0} \label{mod=2} 
\end{align}
for any $g \in H^{-1}(\Omega)$, where $I_h^{ML}$ is the usual interpolant of the conforming linear finite element defined in Section \ref{sec5=5}.

\subsection{Stability}

\begin{lem} \label{lem=6}
Let $\{ \mathbb{T}_h\}$ be a family of conformal meshes with semiregular properties (Assumption \ref{ass1}). Let $T \in \mathbb{T}_h$ be an element of {Condition \ref{cond1}}. We assume that $\Omega$ is convex, and $h \leq 1$. Then, it holds that, for any $\varphi_h \in V_{h0}^M$,
\begin{align}
\displaystyle
| I_h^{ML} \varphi_h |_{H^1(\Omega)} \leq c |\varphi_h|_{H^2(\mathbb{T}_h)}. \label{Merr=1}
\end{align}
\end{lem}

\begin{pf*}
Let  $\varphi_h \in V_{h0}^M$, We set $q:= - \varDelta I_h^{ML} \varphi_h \in H^{-1}(\Omega)$, The dual problem for the fourth-order elliptic problem is to determine $\eta \in H^2_0(\Omega)$ such that:
\begin{align}
\displaystyle
a_1(\varphi , \eta) = \langle q , \varphi \rangle \quad \forall \varphi \in H^2_0(\Omega). \label{Merr=2}
\end{align}
Because $\Omega$ is convex, the solution $\eta$ is in $H^3(\Omega)$, and the following holds: 
\begin{align}
\displaystyle
\| \eta \|_{H^3(\Omega)} \leq c \| q \|_{H^{-1}(\Omega)}. \label{Merr=3}
\end{align}
Let $\mathscr{D}(\Omega)$ be a linear space of infinitely differentiable functions with a compact support on $\Omega$. Using the Gaussian--Green formula, we obtain
\begin{align}
\displaystyle
\langle q , \varphi \rangle = \langle - \varDelta I_h^{ML} \varphi_h , \varphi \rangle
= \int_{\Omega} \nabla I_h^{ML} \varphi_h \cdot \nabla \varphi dx \quad \forall \varphi \in \mathscr{D}(\Omega). \label{Merr=4}
\end{align}
Because $\mathscr{D}(\Omega)$ is dense in $H_0^1(\Omega)$, Equation \eqref{Merr=4} is valid for any $\varphi \in H_0^1(\Omega)$. Then,
\begin{align}
\displaystyle
\| q \|_{H^{-1}(\Omega)}
&= \sup_{0 \neq \varphi \in H_0^1(\Omega)} \frac{\langle q , \varphi \rangle}{| \varphi |_{H^1(\Omega)}} 
= \sup_{0 \neq \varphi \in H_0^1(\Omega)} \frac{\int_{\Omega} \nabla I_h^{ML} \varphi_h \cdot \nabla \varphi dx}{| \varphi |_{H^1(\Omega)}} \notag \\
&\leq | I_h^{ML} \varphi_h |_{H^1(\Omega)}. \label{Merr=5}
\end{align}
Moreover, using \eqref{Merr=2} and the Gauss--Green formula, it holds that, for any $\varphi \in \mathscr{D}(\Omega)$,
\begin{align}
\displaystyle
\langle q , \varphi \rangle
&= \sum_{i=1}^2 \int_{\Omega} \nabla \frac{\partial \eta}{\partial x_i } \cdot \nabla \frac{\partial \varphi}{\partial x_i} dx
=  -  \sum_{i=1}^2 \int_{\Omega} \varDelta \frac{\partial \eta}{\partial x_i } \frac{\partial \varphi}{\partial x_i} dx.  \label{Merr=6}
\end{align}
Equality \eqref{Merr=6} is valid for any $\varphi \in H_0^1(\Omega)$, and by substituting $I_h^{ML} \varphi_h$ for $\varphi$ in \eqref{Merr=4} and \eqref{Merr=6} and using  \eqref{Mh=5}, we obtain for any $v_h \in V_h^{RT^0}$.
\begin{align}
\displaystyle
| I_h^{ML} \varphi_h |^2_{H^1(\Omega)}
&= \langle q , I_h^{ML} \varphi_h \rangle = - \sum_{i=1}^2 \int_{\Omega} \varDelta \frac{\partial \eta}{\partial x_i } \frac{\partial ( I_h^{ML} \varphi_h)}{\partial x_i} dx \notag \\
&\hspace{-2cm} = -  \sum_{T \in \mathbb{T}_h} \sum_{i=1}^2 \int_{T} \varDelta \frac{\partial \eta}{\partial x_i } \frac{\partial \varphi_h}{\partial x_i} dx + \sum_{T \in \mathbb{T}_h} \sum_{i=1}^2 \int_{T} \varDelta \frac{\partial \eta}{\partial x_i } \frac{ \partial (\varphi_h - I_h^{ML} \varphi_h)}{\partial x_i} dx \notag \\
&\hspace{-2cm} =  \sum_{T \in \mathbb{T}_h} \sum_{i=1}^2 \int_{T} \left( v_h - \nabla \frac{\partial \eta}{\partial x_i } \right) \cdot \nabla \frac{\partial \varphi_h}{\partial x_i} dx  \notag\\
&\hspace{-2cm} \quad + \sum_{T \in \mathbb{T}_h} \sum_{i=1}^2 \int_{T} \left( \div v_h - \varDelta \frac{\partial \eta}{\partial x_i } \right) \frac{\partial \varphi_h}{\partial x_i} dx \notag \\
&\hspace{-2cm} \quad + \sum_{T \in \mathbb{T}_h} \sum_{i=1}^2 \int_{T} \varDelta \frac{\partial \eta}{\partial x_i } \frac{ \partial (\varphi_h - I_h^{ML} \varphi_h)}{\partial x_i} dx + \sum_{T \in \mathbb{T}_h} \sum_{i=1}^2 \int_{T}  \nabla \frac{\partial \eta}{\partial x_i } \cdot \nabla \frac{\partial \varphi_h}{\partial x_i} dx \notag \\
&\hspace{-2cm} = I_1 + I_2 + I_3 + I_4.  \label{Merr=7}
\end{align}
We set $v_h := I_h^{RT^0} \nabla \frac{\partial \eta}{\partial x_i }$.  Using \eqref{RT5}, \eqref{Merr=3}, \eqref{Merr=5}, Assumption \ref{ass1}, and the Jenssen inequality, we can obtain an estimate of the term $I_1$ such that
\begin{align}
\displaystyle
|I_1|
&\leq c \sum_{T \in \mathbb{T}_h} \sum_{i=1}^2 \left \|   I_{T}^{RT^0} \nabla \frac{\partial \eta}{\partial x_i } - \nabla \frac{\partial \eta}{\partial x_i } \right \|_{L^2(T)^2} |\varphi_h|_{H^2(T)} \notag \\
&\leq c \sum_{T \in \mathbb{T}_h} \sum_{i=1}^2 \left( \sum_{j=1}^2 h_j \left \|  \frac{\partial}{\partial r_j} \nabla \frac{\partial \eta}{\partial x_i } \right \|_{L^2(T)^2} + h_T \left \| \varDelta \frac{\partial \eta}{\partial x_i} \right \|_{L^{2}({T})} \right) |\varphi_h|_{H^2(T)} \notag \\
&\leq c  \sum_{T \in \mathbb{T}_h} \left( \sum_{j=1}^2 h_j \left| \frac{\partial \eta}{\partial r_j} \right|_{H^2(T)} + h_T |\eta|_{H^3(T)} \right) |\varphi_h|_{H^2(T)} \notag \\
&\leq c \sum_{T \in \mathbb{T}_h} \sum_{j=1}^2 h_j \left| \frac{\partial \eta}{\partial r_j} \right|_{H^2(T)}  |\varphi_h|_{H^2(\mathbb{T}_h)} + h |\eta|_{H^3(T)}  |\varphi_h|_{H^2(\mathbb{T}_h)} \notag \\
&\leq c h |\eta|_{H^3(\Omega)}  |\varphi_h|_{H^2(\mathbb{T}_h)} \leq c h | I_h^{ML} \varphi_h |_{H^1(\Omega)}   |\varphi_h|_{H^2(\mathbb{T}_h)}. \label{Merr=8}
\end{align}
Using \eqref{L2ortho}, \eqref{Merr=3}, and \eqref{Merr=5}, and the stability of the $L^2$projection, we estimate $I_2$ as
\begin{align}
\displaystyle
|I_2|
&= \left|  \sum_{T \in \mathbb{T}_h} \sum_{i=1}^2 \int_{T} \left(  \Pi_{T}^0 \varDelta \frac{\partial \eta}{\partial x_i } - \varDelta \frac{\partial \eta}{\partial x_i} \right) \left( \frac{\partial \varphi_h}{\partial x_i}  - \Pi_{T}^0 \frac{\partial \varphi_h}{\partial x_i} \right) dx \right| \notag \\
&\leq \sum_{T \in \mathbb{T}_h} \sum_{i=1}^2 \left \|  \Pi_{T}^0 \varDelta \frac{\partial \eta}{\partial x_i } - \varDelta \frac{\partial \eta}{\partial x_i } \right \|_{L^2(T)} \left \| \frac{\partial \varphi_h}{\partial x_i} - \Pi_{T}^0 \frac{\partial \varphi_h}{\partial x_i} \right \|_{L^2(T)} \notag \\
&\leq c h |\eta|_{H^3(\Omega)}  |\varphi_h|_{H^2(\mathbb{T}_h)} \leq c h | I_h^{ML} \varphi_h |_{H^1(\Omega)}   |\varphi_h|_{H^2(\mathbb{T}_h)}, \label{Merr=9}
\end{align}
{
where we used the fact that
\begin{align*}
\displaystyle
\div  I_h^{RT^0} \nabla \frac{\partial \eta}{\partial x_i } = \Pi_h^0 \div \nabla \frac{\partial \eta}{\partial x_i } = \Pi_h^0 \varDelta \frac{\partial \eta}{\partial x_i }.
\end{align*}
}
Using the estimates \eqref{error=22} and \eqref{Merr=3}, the term $I_3$ is estimated as follows:
\begin{align}
\displaystyle
|I_3|
&\leq c \sum_{T \in \mathbb{T}_h} |\eta|_{H^3(T)} \sum_{i=1}^2 \left \|  \frac{\partial (I_h^{ML} \varphi_h - \varphi_h)}{\partial x_i} \right\|_{L^2(T)} \notag \\
&\leq c |\eta|_{H^3(\Omega)} |I_h^{ML} \varphi_h - \varphi_h |_{H^1(\mathbb{T}_h)} \notag \\
&\leq c h |\eta|_{H^3(\Omega)} |\varphi_h|_{H^2(\mathbb{T}_h)} \leq c h | I_h^{ML} \varphi_h |_{H^1(\Omega)}   |\varphi_h|_{H^2(\mathbb{T}_h)}. \label{Merr=10}
\end{align}
By substituting $\eta$ for $\varphi$ in \eqref{Merr=4} and \eqref{Merr=5}, we obtain:
\begin{align}
\displaystyle
|\eta|^2_{H^2(\Omega)}
&\leq a(\eta,\eta) = \langle q , \eta \rangle \leq c \| q \|_{H^{-1}(\Omega)} |\eta|_{H^1(\Omega)} \leq c  | I_h^{ML} \varphi_h |_{H^1(\Omega)} |\eta|_{H^1(\Omega)}. \label{Merr=11}
\end{align}
Furthermore, using the Gauss--Green formula and Poincar\'e inequality, we obtain
\begin{align}
\displaystyle
|\eta|^2_{H^1(\Omega)}
&\leq \| \varDelta \eta \|_{L^2(\Omega)} \| \eta \|_{L^2(\Omega)} \leq c |\eta|_{H^2(\Omega)} |\eta|_{H^1(\Omega)}. \label{Merr=12}
\end{align}
From \eqref{Merr=11} and \eqref{Merr=12}, we have that
\begin{align}
\displaystyle
|\eta|_{H^2(\Omega)} \leq c | I_h^{ML} \varphi_h |_{H^1(\Omega)}. \label{Merr=13}
\end{align}
Using \eqref{Merr=13}, we estimate the term $I_4$ as
\begin{align}
\displaystyle
|I_4|
&\leq c |\eta|_{H^2(\Omega)} |\varphi_h|_{H^2(\mathbb{T}_h)} \leq c | I_h^{ML} \varphi_h |_{H^1(\Omega)}  |\varphi_h|_{H^2(\mathbb{T}_h)}. \label{Merr=14}
\end{align}
If $h \leq 1$, combining  \eqref{Merr=7}, \eqref{Merr=8}, \eqref{Merr=9}, \eqref{Merr=10}, and \eqref{Merr=14} yields the target inequality \eqref{Merr=1}.
\qed
\end{pf*}

\begin{thr} \label{sta=th9}
Let $\{ \mathbb{T}_h\}$ be a family of conformal meshes with semiregular properties (Assumption \ref{ass1}). Let $T \in \mathbb{T}_h$ be an element of {Condition \ref{cond1}}. We assume that $\Omega$ is convex, and $h \leq 1$. For any $g \in H^{-1}(\Omega)$, let $\psi_h^* \in V_{h0}^M$ be the solution to \eqref{mod=2}. Subsequently, it holds that
\begin{align}
\displaystyle
|\psi_h^*|_{H^2(\mathbb{T}_h)} \leq c \| g \|_{H^{-1}(\Omega)}. \label{Merr=15}
\end{align}
\end{thr}

\begin{pf*}
By substituting $\psi_h^*$ for $\varphi_h$ in \eqref{mod=2}, and using \eqref{Merr=3}, we determine that the following holds:
\begin{align*}
\displaystyle
|\psi_h^*|_{H^2(\mathbb{T}_h)}^2
&\leq c \| g \|_{H^{-1}(\Omega)} | I_h^{ML} \psi_h^*|_{H^1(\Omega)} \leq c  \| g \|_{H^{-1}(\Omega)} | \psi_h^*|_{H^2(\mathbb{T}_h)},
\end{align*}
which leads to the target inequality \eqref{Merr=15}.
\qed
\end{pf*}

\begin{note}
In Theorem \ref{sta=th9}, we have derived the stability estimate by imposing that $\Omega$ is convex. Because we used the regularity of the solution of the dual problem, the convexity can be not violated in our method.
\end{note}

\subsection{Error analysis} \label{error=analysis64}
The starting point for error analysis of the modified Morley FEM is the following inequality:
\begin{lem}
Let $\psi  \in H^3(\Omega) \cap H^2_0(\Omega)$ be the solution to the fourth-order {elliptic} problem in \eqref{cont=2} with $\varDelta^2 \psi \in H^{-1}(\Omega)$. Let $\psi_h^* \in V_{h0}^M$ be an approximate solution to the problem \eqref{mod=2}. It holds that
\begin{align}
\displaystyle
| \psi - \psi_h^* |_{H^2(\mathbb{T}_h)} 
&\leq c \inf_{v_h \in V_{h0}^M} | \psi -v_h |_{H^2(\mathbb{T}_h)} \notag\\
&\quad + c \sup_{\varphi_h \in V_{h0}^M} \frac{\left|a_{h1}(\psi,\varphi_h) - \int_{\Omega} \varDelta^2 \psi I_h^{ML} \varphi_h dx \right |}{| \varphi_h|_{H^2(\mathbb{T}_h)}}. \label{Merr=16}
\end{align}

\end{lem}

\begin{pf*}
The proof is standard.
\qed
\end{pf*}

The first term on the R.H.S. of inequality \eqref{Merr=16} is estimated as follows: Using the Morley interpolation error estimate \eqref{Mor24} for any $\psi \in H^3(\Omega)$,
\begin{align}
\displaystyle
\inf_{v_h \in V_{h0}^M} | \psi -v_h |_{H^2(\mathbb{T}_h)}
&\leq | \psi -I_h^{M} \psi |_{H^2(\mathbb{T}_h)} \leq c \sum_{T \in \mathbb{T}_h} \sum_{i=1}^2  h_i \left | \frac{\partial \psi}{\partial r_i} \right |_{H^{2}(T)}. \label{Merr=17}
\end{align}

Here, we present the error estimate for the consistency term. It should be noted that $\psi  \in H^3(\Omega) \cap H^2_0(\Omega)$ with $\varDelta^2 \psi \in L^2(\Omega)$.

\begin{lem} \label{lem=new6}
Let $\{ \mathbb{T}_h\}$ be a family of conformal meshes with semiregular properties (Assumption \ref{ass1}). Let $T \in \mathbb{T}_h$ be an element that satisfies Condition \ref{cond1}. Subsequently, there exists a positive constant $C$ independent of $h$, such that for any $\psi \in H^3(\Omega) \cap H^2_0(\Omega)$ with $\varDelta^2 \psi \in L^2(\Omega)$ and $\varphi_h \in V_{h0}^M$,
\begin{align}
\displaystyle
&\left |a_{h1}(\psi,\varphi_h) -\int_{\Omega} \varDelta^2 \psi  I_h^{ML} \varphi_h dx \right| \notag \\
&\leq c \left(  \sum_{T \in \mathbb{T}_h} \sum_{i=1}^2 h_i \left| \frac{\partial \psi}{\partial r_i} \right|_{H^2(T)} \right)  |\varphi_h|_{H^2(\mathbb{T}_h)} + c h \| \nabla \varDelta \psi \|_{L^2(\Omega)^2} |\varphi_h|_{H^2(\mathbb{T}_h)}. \label{Merr=18}
\end{align}
\end{lem}

\begin{pf*}
Let $\psi \in H^3(\Omega) \cap H^2_0(\Omega)$ with $\varDelta^2 \psi \in L^2(\Omega)$ and $\varphi_h \in V_{h0}^M $. Applying the Gaussian--Green formula and \eqref{Mh=5} yields, for any $v_{h,i} \in V_h^{RT^0}$, $i=1,2$,
\begin{align}
\displaystyle
&a_{h1}(\psi,\varphi_h) - \int_{\Omega} \varDelta^2 \psi  I_h^{ML} \varphi_h  dx \notag\\
&=  \sum_{T \in \mathbb{T}_h} \sum_{i=1}^2 \int_{T} \left( \nabla \frac{\partial \psi}{\partial x_i } \cdot \nabla \frac{\partial \varphi_h}{\partial x_i} + \varDelta \frac{\partial \psi}{\partial x_i} \frac{\partial \varphi_h}{\partial x_i} + \varDelta \frac{\partial \psi}{\partial x_i} \frac{\partial (I_h^{ML} \varphi_h - \varphi_h)}{\partial x_i} \right) dx \notag \\
&= \sum_{T \in \mathbb{T}_h} \sum_{i=1}^2 \int_{T} \left( \nabla \frac{\partial \psi}{\partial x_i } - v_{h,i} \right) \cdot \nabla \frac{\partial \varphi_h}{\partial x_i} dx \notag \\
&+  \sum_{T \in \mathbb{T}_h} \sum_{i=1}^2 \int_{T} \left( \varDelta \frac{\partial \psi}{\partial x_i}  - \div v_{h,i} \right) \frac{\partial \varphi_h}{\partial x_i} dx
+ \sum_{T \in \mathbb{T}_h} \sum_{i=1}^2 \int_{T} \varDelta \frac{\partial \psi}{\partial x_i} \frac{\partial (I_h^{ML} \varphi_h - \varphi_h)}{\partial x_i} dx \notag \\
&=: J_1 + J_2 + J_3. \label{Merr=19}
\end{align}
Here, we use the first
\begin{align*}
\displaystyle
- \int_{\Omega} \varDelta^2 \psi  I_h^{ML} \varphi_h  dx
&= - \int_{\partial \Omega} (n \cdot \nabla) \varDelta \psi I_h^{ML} \varphi_h ds + \int_{\Omega} \nabla  \varDelta \psi \cdot \nabla I_h^{ML} \varphi_h dx \\
&= \int_{\Omega} \nabla  \varDelta \psi \cdot \nabla I_h^{ML} \varphi_h dx,
\end{align*}
because $ I_h^{ML} \varphi_h \in H_0^1(\Omega)$,

In setting $v_{h,i} := I_h^{RT^0} \nabla \frac{\partial \psi}{\partial x_i }$, the same method estimates terms $J_1$ and $J_2$ as terms $I_1$ and $I_2$ in Lemma \ref{lem=6}. Subsequently,
\begin{align}
\displaystyle
|J_1|
&\leq c \sum_{T \in \mathbb{T}_h} \sum_{j=1}^2 h_j \left| \frac{\partial \psi}{\partial r_j} \right|_{H^2(T)}  |\varphi_h|_{H^2(\mathbb{T}_h)} + c h \| \nabla \varDelta \psi \|_{L^2(\Omega)^2}  |\varphi_h|_{H^2(\mathbb{T}_h)}, \label{Merr=20} \\
|J_2|
&\leq c h \| \nabla \varDelta \psi \|_{L^2(\Omega)^2} |\varphi_h|_{H^2(\mathbb{T}_h)}. \label{Merr=21}
\end{align}
Using the estimate \eqref{error=22}, term $J_3$ is estimated using the same method as that for term $I_3$ in lemma \ref{lem=6}:
\begin{align}
\displaystyle
|J_3|
&\leq c h \| \nabla \varDelta \psi \|_{L^2(\Omega)^2} |\varphi_h|_{H^2(\mathbb{T}_h)}. \label{Merr=22}
\end{align}
Using \eqref{Merr=19}, \eqref{Merr=20}, \eqref{Merr=21}, and \eqref{Merr=22}, we can deduce the desired estimate \eqref{Merr=18}.
\qed
\end{pf*}

We {present} an error estimate for a modified Morley FEM for a fourth-order {elliptic} problem.

\begin{thr} \label{Mod=thr6}
Let $\{ \mathbb{T}_h\}$ be a family of conformal meshes with semiregular properties (Assumption \ref{ass1}). We assume all elements $T \in \mathbb{T}_h$ satisfy Condition \ref{cond1}. Let $\psi  \in H^3(\Omega) \cap H^2_0(\Omega)$, where $\varDelta^2 \psi \in L^2(\Omega)$ is the solution of the fourth-order {elliptic} problem in \eqref{cont=2}. Let $\psi_h^* \in V_{h0}^M$ be an approximate solution to the problem \eqref{mod=2}. Subsequently, 
\begin{align}
\displaystyle
| \psi - \psi_h^* |_{H^2(\mathbb{T}_h)} &\leq c \sum_{T \in \mathbb{T}_h} \sum_{i=1}^2  h_i \left | \frac{\partial \psi}{\partial r_i} \right |_{H^{2}(T)} + c h  \| \nabla \varDelta \psi \|_{L^2(\Omega)^2}. \label{Merr=23}
\end{align}
\end{thr}

\begin{pf*}
By combining \eqref{Merr=16}, \eqref{Merr=17}, and \eqref{Merr=18}, we obtain target estimate \eqref{Merr=23}.
\qed
\end{pf*}

\subsection{Usual Morley finite element method}
From Theorem \ref{Mod=thr6}, we can deduce the error estimate for a typical Moley finite element problem \eqref{mod=1}.

\begin{coro} \label{coro=1}
Let $\{ \mathbb{T}_h\}$ be a family of conformal meshes with semiregular properties (Assumption \ref{ass1}). We assume all elements $T \in \mathbb{T}_h$ satisfy Condition \ref{cond1}. Let $\psi  \in H^3(\Omega) \cap H^2_0(\Omega)$, where $\varDelta^2 \psi \in L^2(\Omega)$ is the solution of the fourth-order {elliptic} problem in \eqref{cont=2}. Let $\psi_h \in V_{h0}^M$ be an approximate solution to the problem \eqref{mod=1}. Then, 
\begin{align}
\displaystyle
| \psi - \psi_h |_{H^2(\mathbb{T}_h)} &\leq c \sum_{T \in \mathbb{T}_h} \sum_{i=1}^2  h_i \left | \frac{\partial \psi}{\partial r_i} \right |_{H^{2}(T)} + c h  \| \nabla \varDelta \psi \|_{L^2(\Omega)^2} + c \nu^{-1} h^2 \| g \|_{L^{2}(\Omega)}. \label{Merr=24}
\end{align}
\end{coro}

\begin{pf*}
The subtraction of \eqref{mod=1} from \eqref{mod=2} with $\varphi_h := \psi_h^* - \psi_h$ and using \eqref{error=23} yields
\begin{align*}
\displaystyle
\nu |\psi_h^* - \psi_h|^2_{H^2(\mathbb{T}_h)}
&\leq \nu \| g \|_{L^2(\Omega)} \| I_h^{ML} ( \psi_h^* - \psi_h) - ( \psi_h^* - \psi_h) \|_{L^2(\Omega)} \\
&\leq c \nu h^2 \| g \|_{L^2(\Omega)} |\psi_h^* - \psi_h|_{H^2(\mathbb{T}_h)},
\end{align*}
which leads to
\begin{align*}
\displaystyle
| \psi - \psi_h |_{H^2(\mathbb{T}_h)} 
&\leq | \psi - \psi_h^* |_{H^2(\mathbb{T}_h)} + | \psi_h^* - \psi_h |_{H^2(\mathbb{T}_h)} \\
&\leq c \sum_{T \in \mathbb{T}_h} \sum_{i=1}^2  h_i \left | \frac{\partial \psi}{\partial r_i} \right |_{H^{2}(T)} + c h \| \nabla \varDelta \psi \|_{L^2(\Omega)^2} + c \nu^{-1} h^2 \| g \|_{L^{2}(\Omega)}.
\end{align*}
\qed
\end{pf*}

\subsection{Three-dimensional case}
When $d=3$, we {consider} the error estimate for the consistency term of a typical Morley FEM.  \\
\indent
The starting point for the error analysis of the classical Morley FEM is as follows. Let $\psi  \in H^3(\Omega) \cap H^2_0(\Omega)$ be the solution to the fourth-order {elliptic} problem \eqref{cont=2} with $\varDelta^2 \psi \in L^2(\Omega)$. Let $\psi_h \in V_{h0}^M$ be an approximate solution to the problem \eqref{mod=1}. It holds that
\begin{align}
\displaystyle
| \psi - \psi_h |_{H^2(\mathbb{T}_h)} 
&\leq c \inf_{v_h \in V_{h0}^M} | \psi -v_h |_{H^2(\mathbb{T}_h)} \notag\\
&\quad + c \sup_{\varphi_h \in V_{h0}^M} \frac{\left|a_{h1}(\psi,\varphi_h) - \int_{\Omega} \varDelta^2 \psi \varphi_h dx \right |}{| \varphi_h|_{H^2(\mathbb{T}_h)}}. \label{newMerr=16}
\end{align}
We {use} the RT interpolation error estimates for the error analysis of the consistency term (see Section \ref{error=analysis64}). However, as stated in Note \ref{newnote5}, we cannot obtain anisotropic RT interpolation error estimates for (Type \roman{stwo}) in Section \ref{sec322}. We use only case  (Type \roman{sone}). Note that standard error analysis for (Type \roman{stwo}) can obtain the error estimates on isotropic meshes. \\
\indent
Let $\{ \mathbb{T}_h\}$ be a family of conformal meshes with semiregular properties (Assumption \ref{ass1}). Let $T \in \mathbb{T}_h$ be an element satisfying Condition \ref{cond2} which satisfy (type \roman{sone}) in Section \ref{sec322}. Let $w \in H^1_0(\Omega)$. We then have
\begin{align}
\displaystyle
a_{h1}(\psi,\varphi_h) - \int_{\Omega} \varDelta^2 \psi \varphi_h dx
= a_{h1}(\psi,\varphi_h) - \int_{\Omega} \varDelta^2 \psi  w dx + \int_{\Omega} \varDelta^2 \psi  (w -  \varphi_h) dx. \label{newerror=630}
\end{align}
For the first and second terms on the R.H.S. of \eqref{newerror=630}, the Gauss--Green formula and \eqref{Mh=5} yield:
\begin{align}
\displaystyle
&a_{h1}(\psi,\varphi_h) - \int_{\Omega} \varDelta^2 \psi  w dx \notag\\
&\quad = \sum_{T \in \mathbb{T}_h} \sum_{i=1}^3 \int_{T} \left( \nabla \frac{\partial \psi}{\partial x_i } \cdot \nabla \frac{\partial \varphi_h}{\partial x_i} + \varDelta \frac{\partial \psi}{\partial x_i} \frac{\partial \varphi_h}{\partial x_i} + \varDelta \frac{\partial \psi}{\partial x_i} \frac{\partial (w - \varphi_h)}{\partial x_i} \right) dx \notag \\
&\quad= \sum_{T \in \mathbb{T}_h} \sum_{i=1}^3 \int_{T} \left( \nabla \frac{\partial \psi}{\partial x_i } - I_h^{RT^0} \nabla \frac{\partial \psi}{\partial x_i } \right) \cdot \nabla \frac{\partial \varphi_h}{\partial x_i} dx \notag \\
&\quad\quad +  \sum_{T \in \mathbb{T}_h} \sum_{i=1}^3 \int_{T} \left( \varDelta \frac{\partial \psi}{\partial x_i}  - \div I_h^{RT^0} \nabla \frac{\partial \psi}{\partial x_i } \right) \frac{\partial \varphi_h}{\partial x_i} dx \notag \\
&\quad\quad + \sum_{T \in \mathbb{T}_h} \sum_{i=1}^3 \int_{T} \varDelta \frac{\partial \psi}{\partial x_i} \frac{\partial (w - \varphi_h)}{\partial x_i} dx. \label{newerror=631}
\end{align}
{The first and second terms on the R.H.S. of \eqref{newerror=631} can be estimated in the same way as the terms $J_1$ and $J_2$ in \eqref{Merr=19}.} To complete the consistency error estimate, we must estimate the third term of the R.H.S. of \eqref{newerror=630} and the third term of the R.H.S. of \eqref{newerror=631}; that is, we show that there exists $w \in H_0^1(\Omega)$ such that
\begin{align}
\displaystyle
\left \| \frac{\partial (w - \varphi_h)}{\partial x_i} \right \|_{L^2(T)} &\leq c h |\varphi_h|_{H^2(T)}, \quad i \in \{1,2,3 \}, \label{error=17} \\
\| w -  \varphi_h \|_{L^2(T)}
 &\leq c h^2  |\varphi_h|_{H^2(T)}. \label{error=18}
\end{align}
However, achieving \eqref{error=17} and \eqref{error=18} for anisotropic meshes may be difficult. Meanwhile, we can deduce the error estimate of the usual Morley finite element method on isotropic meshes, see \cite[Lemma 6]{WanJin06}.

{
\begin{rem}
Several problems can lead to fourth-order (partial) differential equations in three dimensions. For example, a fourth-order term appears in the following real-world problem: The Cahn--Hilliard equation \cite{CahHil58}, which describes the phase separation phenomenon between two metals, and the bending of thin plates, which is described by the biharmonic equation.
\end{rem}
}

\section{Application to stream function formulation}

\subsection{Continuous problem}
Let $\Omega \subset \mathbb{R}^d$ and $d \in \{ 2,3\}$ be {a bounded polyhedral domain}. The (scaled) Stokes equation is as follows: Determine $(u,p): \Omega \to \mathbb{R}^d \times \mathbb{R}$ such that
\begin{align}
\displaystyle
- \nu \varDelta u + \nabla p = f \quad \text{in $\Omega$}, \quad \div u  = 0 \quad \text{in $\Omega$}, \quad u = 0 \quad \text{on $\partial \Omega$}, \label{cont=3}
\end{align}
where $\nu$ is a nonnegative parameter and $f:\Omega \to \mathbb{R}^d$ is a given function. The variational formulation for the standard Stokes equation is as follows: For any $f \in L^2(\Omega)^d$, determine $(u,p) \in V:= H^1_0(\Omega)^d \times Q:= L^2_0(\Omega):=  \left \{ q \in L^2(\Omega); \ \int_{\Omega} q dx = 0 \right \}$ such that
\begin{subequations} \label{cont=4}
\begin{align}
\displaystyle
\nu a_2(u,v) + b(v , p) &= \int_{\Omega} f \cdot v dx \quad \forall v \in V, \label{cont=4a} \\
b(u , q) &= 0 \quad \forall q \in Q, \label{cont=4b}
\end{align}
\end{subequations}
where the bilinear forms $a_2: H^{1}(\Omega)^d \times H^1(\Omega)^d \to \mathbb{R}$ and $b:H^1(\Omega)^d \times L^2(\Omega) \to \mathbb{R}$ are
\begin{align*}
\displaystyle
&a_2(u , v) :=  \sum_{i=1}^d \int_{\Omega} \nabla u_i \cdot \nabla v_i dx, \quad b(v,q) := - \int_{\Omega} \div v q dx.
\end{align*}
for any $u,v \in H^1(\Omega)^d$ or $q \in L^2(\Omega)$. Using the space of weakly divergence-free functions, we obtain
\begin{align*}
\displaystyle
V_{\div} := \{ v \in V  ; \ b(v,q) = 0 \ \forall q \in Q \} \subset V,
\end{align*}
The problem associated with \eqref{cont=4} is as follows. We determine $u \in V_{\div}$ such that
\begin{align}
\displaystyle
\nu a_2(u,v) &= \int_{\Omega} f \cdot v dx \quad \forall v \in V_{\div}. \label{cont=5}
\end{align}

\begin{thr}
If $d=2$ and $\Omega$ {is} convex, then the solution $(u,p)$ to the Stokes problem belongs to $u \in H^2(\Omega)^2$ and $p \in H^1(\Omega)$.
\end{thr}

\begin{pf*}
Reference \cite{KelOsb76} provides the proof.
\qed
\end{pf*}

Let $d=2$ and assume that $\Omega$ is connected. As $u \in H^1_0(\Omega)^2$ is divergence-free, it holds that $u = \curl \psi$ for a unique stream function $\psi \in H^2_0(\Omega)$; see \cite[Section \Roman{lone}.3.1]{GirRav86}. By setting $u := \curl \psi$, Stokes problem {\eqref{cont=3}} is reduced to 
\begin{align}
\displaystyle
\nu \varDelta^2 \psi &= g:= \rot f = \frac{\partial f_2}{\partial x_1} - \frac{\partial f_1}{\partial x_2} \quad \text{in $\Omega$}, \quad \psi = \frac{\partial \psi}{\partial n} = 0 \quad \text{on $\partial \Omega$}, \label{stream}
\end{align}
for any $f = (f_1,f_2)^T \in L^2(\Omega)^2$. This function is given in $H^{-1}(\Omega)$; see \cite[Section \Roman{lone}.5.2]{GirRav86} and \cite[Section 7.3.3]{Gri11}. Moreover, by setting $u := \curl \psi$ and $v := \curl \varphi$, the Stokes problem {\eqref{cont=5}} can be reduced to 
\begin{align}
\displaystyle
\nu a_1(\psi,\varphi) &= \int_{\Omega} f \cdot \curl \varphi dx \quad \forall \varphi \in H^2_0(\Omega). \label{stream2}
\end{align}

\subsection{Modified finite element method}
The CR finite element space is defined as follows:
\begin{align*}
\displaystyle
V_{h0}^{CR} &:=  \biggl \{ \varphi_h \in P_{dc,h}^1: \  \int_F [\![ \varphi_h ]\!] ds = 0 \ \forall F \in \mathcal{F}_h \biggr \}.
\end{align*}
We then define the Stokes pair $(V_h,Q_h)$ as
\begin{align*}
\displaystyle
V_h := (V_{h0}^{CR})^d, \quad Q_h := P_{dc,h}^0 \cap L^{2}_0(\Omega).
\end{align*}
Furthermore, we define  $a_{h2}: (H^1(\Omega)^d + V_{h}) \times H^1 (\Omega)^d + V_{h}) \to \mathbb{R}$ and $b_h: (V+V_{h}) \times Q_h \to \mathbb{R}$, which are the discrete counterparts of the bilinear forms $a_{2}$ and $b$ as follows:
\begin{align*}
\displaystyle
a_{h2}(u_h,v_h) &:= \sum_{i=1}^d \sum_{T \in \mathbb{T}_h} \int_{T} \nabla u_{h,i} \cdot \nabla v_{h,i} dx, \\
b_h(v_h , q_h) &:= - \sum_{T \in \mathbb{T}_h} \int_{T} \div v_h q_h dx.
\end{align*}
The modified CR finite element approximation problem for the Stokes equation proposed in \cite{Lin14} is as follows. Determine $(u_h^*,p_h^*) \in V_{h} \times Q_h$ such that
\begin{subequations} \label{mod=8}
\begin{align}
\displaystyle
\nu a_{h2}(u_h^*,v_h) + b_h(v_h , p_h^*) &= \int_{\Omega} f \cdot \mathcal{I}_{h}^{RT^0} v_h dx\quad \forall v_h \in V_{h}, \label{mod=8a} \\
b_h(u_h^* , q_h) &= 0 \quad \forall q_h \in Q_h, \label{mod=8b}
\end{align}
\end{subequations}
where the lifting operator $\mathcal{I}_{h}^{RT^0}$ is defined in Section \ref{sec=RT}. We define a discrete weakly divergence-free subspace as
\begin{align*}
\displaystyle
V_{h,\div} := \{ v_h \in V_h: \ b_h(v_h,q_h) = 0 \ \forall q_h \in Q_h \}.
\end{align*}
Because $V_{h,\div} \not\subset V_{\div}$, {the space $V_{h,\div}$} is nonconforming in {the space $V_{\div}$}. Subsequently, the problem associated with \eqref{mod=8} is to determine $u_h^* \in V_{h,\div}$ such that:
\begin{align}
\displaystyle
\nu a_{h2}(u_h^*,v_h) = \int_{\Omega} f \cdot \mathcal{I}_{h}^{RT^0} v_h dx \quad \forall v_h \in V_{h,\div}. \label{mod=9}
\end{align}

\subsection{Error analysis}
Let $d=2$ and assume that $\Omega$ is connected. For any $\varphi_h \in V_{h0}^M$, we define the broken curl as
\begin{align*}
\displaystyle
(\curlh \varphi_h)|_T = \curl \varphi_h|_T = \left( \frac{\partial \varphi_h}{\partial x_2}, -  \frac{\partial \varphi_h}{\partial x_1} \right){^{\top}} \Biggl |_T \quad \forall T \in \mathbb{T}_h.
\end{align*}
The following theorem is known.

\begin{thr}
The following holds:
\begin{align*}
\displaystyle
V_{h,\div} = \curlh V_{h0}^M.
\end{align*}

\end{thr}

\begin{pf*}
Reference \cite[Theorem 4.1]{FalMor90} provides the proof.
\qed
\end{pf*}

The modified Morley FEM from the stream function formulation is as follows: For any $\psi_h^{**}, \varphi_h \in V_{h0}^M$, substituting $u_h^*:= \curlh \psi_h^{**}$ and $v_h:= \curlh \varphi_h$ in the problem \eqref{mod=9} yields the following problem: For any $f = (f_1,f_2){^{\top}} \in L^2(\Omega)^2$, determine $\psi_h^{**} \in V_{h0}^M$ such that
\begin{align}
\displaystyle
\nu a_{h1}(\psi_h^{**} , \varphi_h) =  \int_{\Omega} f \cdot I_{h}^{RT^0} (\curlh \varphi_h) dx \quad \forall \varphi_h \in V_{h0}^M. \label{mod=10}
\end{align}

We {consider} a typical CR finite element method for the Stokes equation. 
\begin{align}
\displaystyle
\nu a_{h1}(\psi_h^{**} , \varphi_h) =  \int_{\Omega} f \cdot \curlh \varphi_h dx \quad \forall \varphi_h \in V_{h0}^M. \label{mod=10b}
\end{align}
For a sufficiently smooth function $f$, the R.H.S. of \eqref{mod=10b} cannot be replaced with
\begin{align*}
\displaystyle
\sum_{T \in \mathbb{T}_h} \int_T \left( \frac{\partial f_2}{\partial x_1} -  \frac{\partial f_1}{\partial x_2} \right) \varphi_h dx.
\end{align*}
Because we have
\begin{align}
\displaystyle
\sum_{T \in \mathbb{T}_h} \int_T  f \cdot \curlh \varphi_h dx
&= \sum_{T \in \mathbb{T}_h} \int_T \left( \frac{\partial f_2}{\partial x_1} -  \frac{\partial f_1}{\partial x_2} \right) \varphi_h dx \nonumber \\
&\quad + \sum_{T \in \mathbb{T}_h} \int_{\partial T} (f_1 n_T^{(2)} - f_2 n_T^{(1)}) \varphi_h ds, \label{mod=10c}
\end{align}
where $n_T := (n_T^{(1)} , n_T^{(2)}){^{\top}}$ and the second term on the R.H.S. of \eqref{mod=10c} does not generally vanish.

The following inequality is the starting point for the error analysis of the modified Morley FEM \eqref{mod=10}:
\begin{lem}
Let $\psi  \in H^3(\Omega) \cap H^2_0(\Omega)$ be the solution to the fourth-order {elliptic} problem \eqref{stream2}. Let $\psi_h^{**} \in V_{h0}^M$ be an approximate solution to \eqref{mod=10}. Then, it holds that
\begin{align*}
\displaystyle
| \psi - \psi_h^{**} |_{H^2(\mathbb{T}_h)} 
&\leq c \inf_{v_h \in V_{h0}^M} | \psi -v_h |_{H^2(\mathbb{T}_h)} \notag\\
&\quad + c \nu^{-1} \sup_{\varphi_h \in V_{h0}^M} \frac{\left|\nu a_{h1}(\psi,\varphi_h) - \int_{\Omega} f \cdot I_{h}^{RT^0} (\curlh \varphi_h) dx \right |}{| \varphi_h|_{H^2(\mathbb{T}_h)}}.
\end{align*}

\end{lem}

\begin{pf*}
The proof is standard.
\qed
\end{pf*}

Let $(u,p)$ be the solution to the problem \eqref{cont=3} with $f \in L^2(\Omega)^2$. Because
\begin{align*}
\displaystyle
\int_{\Omega} \nabla p \cdot I_{h}^{RT^0} (\curlh \varphi_h) dx = 0, \quad u = \curl \psi, \quad \psi \in H^2_0(\Omega),
\end{align*}
and $\mathcal{I}_{h}^{RT^0} v_h$ is divergence-free on $T$ for any $v_h \in V_{h}$ (see \cite[Lemma 4.134]{Joh16}), we have
\begin{align*}
\displaystyle
&\nu a_{h1}(\psi,\varphi_h) - \int_{\Omega} f \cdot I_{h}^{RT^0} (\curlh \varphi_h) dx \\
&\quad = \nu a_{h1}(\psi,\varphi_h) + \int_{\Omega} \nu \varDelta \curl \psi \cdot I_{h}^{RT^0} (\curlh \varphi_h) dx.
\end{align*}

We present only the error estimates for the consistency term.

\begin{lem} \label{lem8}
Let $\{ \mathbb{T}_h\}$ be a family of conformal meshes with semiregular properties (Assumption \ref{ass1}). Let $T \in \mathbb{T}_h$ be an element that satisfies Condition \ref{cond1}. Subsequently, there exists a positive constant $C$ independent of $h$ such that for any $\psi \in H^3(\Omega) \cap H^2_0(\Omega)$ and $\varphi_h \in V_{h0}^M$,
\begin{align}
\displaystyle
&\left |a_{h1}(\psi,\varphi_h) + \int_{\Omega}  \varDelta \curl \psi \cdot I_{h}^{RT^0} (\curlh \varphi_h) dx \right| \notag \\
&\quad \leq c \left(  \sum_{T \in \mathbb{T}_h} \sum_{j=1}^2  h_j \left | \frac{\partial \psi}{\partial r_j} \right |_{H^{2}(T)}  \right)  |\varphi_h|_{H^2(\mathbb{T}_h)} 
 + c h \left \| \varDelta \curl \psi \right \|_{L^{2}({\Omega})^2} |\varphi_h|_{H^2(\mathbb{T}_h)}. \label{Merr=25}
\end{align}
\end{lem}

\begin{pf*}
Let $\psi \in H^3(\Omega) \cap H^2_0(\Omega)$ and $\varphi_h \in V_{h0}^M $. By using the Gaussian--Green formula and \eqref{Mh=5}, for any $v_{h,i} \in V_h^{RT^0}$ and $i=1,2$,
\begin{align}
\displaystyle
&a_{h1}(\psi,\varphi_h) +\int_{\Omega}  \varDelta \curl \psi \cdot I_{h}^{RT^0} (\curlh \varphi_h) dx \notag \\
&\quad= \sum_{i=1}^2 \sum_{T \in \mathbb{T}_h} \int_{T} \nabla \frac{\partial \psi}{\partial x_i} \cdot \nabla \frac{\partial \varphi_h}{\partial x_i} dx +\int_{\Omega}  \varDelta \curl \psi \cdot \curlh \varphi_h dx \notag\\
&\quad \quad + \int_{\Omega}  \varDelta \curl \psi \cdot \left(  I_{h}^{RT^0} (\curlh \varphi_h) -\curlh \varphi_h \right) dx \notag\\
&\quad= \sum_{T \in \mathbb{T}_h} \sum_{i=1}^2 \int_{T} \left( \nabla \frac{\partial \psi}{\partial x_i } - v_{h,i} \right) \cdot \nabla \frac{\partial \varphi_h}{\partial x_i} dx \notag \\
&\quad \quad +  \sum_{T \in \mathbb{T}_h} \sum_{i=1}^2 \int_{T} \left( \varDelta \frac{\partial \psi}{\partial x_i}  - \div v_{h,i} \right) \frac{\partial \varphi_h}{\partial x_i} dx \notag \\
&\quad \quad + \int_{\Omega}  \varDelta \curl \psi \cdot \left(  I_{h}^{RT^0} (\curlh \varphi_h) -\curlh \varphi_h \right) dx \notag \\
&\quad = K_1 + K_2 + K_3. \label{Merr=26}
\end{align}
By setting $v_{h,i} := I_h^{RT^0} \nabla \frac{\partial \psi}{\partial x_i }$ and using \eqref{jensen}, we can estimate terms $K_1$ and $K_2$ using the same method as for terms $I_1$ and $I_2$ in Lemma \ref{lem=6} and $J_1$ and $J_2$ in Lemma \ref{lem=new6}. 
\begin{align}
\displaystyle
|K_1|
&\leq c \sum_{T \in \mathbb{T}_h} \sum_{i=1}^2 \left \| \nabla \frac{\partial \psi}{\partial x_i } -  I_{T}^{RT^0} \nabla \frac{\partial \psi}{\partial x_i }  \right \|_{L^2(T)^2} |\varphi_h|_{H^2(T)} \notag \\
&\leq c  \sum_{T \in \mathbb{T}_h} \left( \sum_{i,j=1}^2 h_j \left \|  \frac{\partial^2 \curl \psi }{\partial x_i \partial r_j} \right \|_{L^2(T)^2} + h_T \left \| \varDelta \curl \psi \right \|_{L^{2}({T})^2} \right) |\varphi_h|_{H^2(T)} \notag \\
&\leq c  \sum_{T \in \mathbb{T}_h} \left( \sum_{j=1}^2 h_j \left | \frac{\partial \psi}{\partial r_j} \right |_{H^{2}(T)} + h_T \left \| \varDelta \curl \psi \right \|_{L^{2}({T})^2} \right) |\varphi_h|_{H^2(T)} \notag \\
&\leq c \left\{ \left(   \sum_{T \in \mathbb{T}_h} \sum_{j=1}^2 h_j \left | \frac{\partial \psi}{\partial r_j} \right |_{H^{2}(T)}  \right) + h \left \| \varDelta \curl \psi \right \|_{L^{2}({\Omega})^2} \right\}   |\varphi_h|_{H^2(\mathbb{T}_h)}, \label{Merr=27}
\end{align}
and
\begin{align}
\displaystyle
|K_2|
&= \left|  \sum_{T \in \mathbb{T}_h} \sum_{i=1}^2 \int_{T} \left(  \Pi_{T}^0 \varDelta \frac{\partial \psi}{\partial x_i } - \varDelta \frac{\partial \psi}{\partial x_i} \right) \left( \frac{\partial \varphi_h}{\partial x_i}  - \Pi_{T}^0 \frac{\partial \varphi_h}{\partial x_i} \right) dx \right| \notag \\
&\leq \sum_{T \in \mathbb{T}_h} \sum_{i=1}^2 \left \|  \Pi_{T}^0 \varDelta \frac{\partial \psi}{\partial x_i } - \varDelta \frac{\partial \psi}{\partial x_i } \right \|_{L^2(T)} \left \| \frac{\partial \varphi_h}{\partial x_i} - \Pi_{T}^0 \frac{\partial \varphi_h}{\partial x_i} \right \|_{L^2(T)} \notag \\
&\leq c h\|  \varDelta \curl \psi  \|_{L^2(\Omega)^2}  |\varphi_h|_{H^2(\mathbb{T}_h)}. \label{Merr=28}
\end{align}
Using \eqref{RT5}, we estimate $K_3$: 
\begin{align}
\displaystyle
|K_3|
&\leq c h \|  \varDelta \curl \psi  \|_{L^2(\Omega)^2} | \curl \varphi_h |_{H^1(\mathbb{T}_h)} \leq c h\|  \varDelta \curl \psi  \|_{L^2(\Omega)^2} | \varphi_h |_{H^2(\mathbb{T}_h)}. \label{Merr=29}
\end{align}
The target inequality is proved based on \eqref{Merr=26}, \eqref{Merr=27}, \eqref{Merr=28}, and \eqref{Merr=29}.
\qed
\end{pf*}

\subsection{Conclusion}
This study {presented} precise anisotropic interpolation error estimates for nonconforming FEMs. As applications, we {demonstrated} the anisotropic error estimates for the fourth-order elliptic problem \eqref{stream} and stream function formulation \eqref{stream2}. We {denoted} the solutions of the typical Morley, {Arnold--Brezzi} modified Morley, and the modified methods \eqref{mod=10} from the stream function formulation by $\psi_h$, $\psi_h^*$, and $\psi_h^{**}$, respectively. For a convex polygonal domain, the following estimates hold.
\begin{align*}
\displaystyle
| \psi - \psi_h |_{H^2(\mathbb{T}_h)} 
&\leq c \sum_{T \in \mathbb{T}_h} \sum_{i=1}^2  h_i \left | \frac{\partial \psi}{\partial r_i} \right |_{H^{2}(T)} + c h  \| \nabla \varDelta \psi \|_{L^2(\Omega)^2}
+ c \frac{h^2}{\nu} \left \| \rot f \right \|_{L^{2}(\Omega)}, \\
| \psi - \psi_h^* |_{H^2(\mathbb{T}_h)} 
&\leq c \sum_{T \in \mathbb{T}_h} \sum_{i=1}^2  h_i \left | \frac{\partial \psi}{\partial r_i} \right |_{H^{2}(T)} + c h  \| \nabla \varDelta \psi \|_{L^2(\Omega)^2}, \\
| \psi - \psi_h^{**} |_{H^2(\mathbb{T}_h)} 
&\leq c \sum_{T \in \mathbb{T}_h} \sum_{i=1}^2  h_i \left | \frac{\partial \psi}{\partial r_i} \right |_{H^{2}(T)} + c h \left \| \varDelta u \right \|_{L^{2}({\Omega})^2}.
\end{align*}
These are more delicate estimates than those of Falk and Morley \cite{FalMor90}, and they also apply without imposing the shape-regularity condition; however, a maximal angle condition is required. However, we impose the regularity assumption $\psi  \in H^3(\Omega) \cap H^2_0(\Omega)$ with $\varDelta^2 \psi \in L^2(\Omega)$ that is stronger than the original {Arnold--Brezzi} modified Morley method to prove the error estimate of the consistency term. The errors may be larger if the scheme is used without appropriate modifications and appropriate mesh partitions are not used. We aim to extend the proposed approach to numerical comparisons, fourth-order elliptic singular perturbation problems, and Navier--Stokes problems.

\begin{acknowledgements}
We thank the reviewers for their valuable comments and suggestions.
\end{acknowledgements}


\end{document}